\documentclass[12pt]{article}
\usepackage{graphics}
\usepackage{amsmath}
\usepackage{amssymb}
\usepackage{graphicx}
\usepackage{makeidx}
\usepackage{mathrsfs}
\usepackage{amsfonts}
\usepackage[mathscr]{eucal}
\usepackage{amsmath}
\usepackage{mathrsfs}
\usepackage{epstopdf}
\usepackage{color}

\def\e{{\rm e}}
\def\e{\hbox{e}}
\def\d{\hbox {d}}
\def\dt{\hbox{dt}}
\def\ds{\displaystyle}
\def\RR{\vbox {\hbox to 8.9pt {I\hskip-2.1pt R\hfil}}}
\def\NN{{\rm I\hskip-2pt N}}
\def\CC{{\rm C\hskip-4.8pt \vrule height 6pt width 12000sp\hskip 5pt}}
\def\q{\quad}  
\def\cen{\centerline}
\def \rec#1{{1\over{#1}}}
\def\pni{\par\noindent}
\def\vsh{\smallskip}
\def\vs{\medskip}

\def\vsp{\vsh\pni} %% ie. \smallskip + \par
\def\vsn{\vsh\pni}
\def\eg{{\it e.g.}\ }

\def\e{\hbox{e}}
\def\exp{\hbox{exp}}
\def\ds{\displaystyle}

\def\q{\quad}	 

\def\lt{\left} \def\rt{\right}

\def\rec#1{\frac{1}{#1}}

\def\Gs{{\widetilde G(s)}}

\def\L{{\cal L}} %%% Laplace Transform !!!!
\def\F{{\cal F}} %%% Fourier Transform !!!!
\def\Gc{{\cal {G}}_c}	 %% CAUCHY PROBLEM
\def\Gs{{\cal {G}}_s}	 %% SIGNALLING PROBLEM
\def\d{\partial}
  \def\dt{\partial t}
\def\dx{\partial x}     
\def\Ai{\hbox{Ai}}
\def\MM{\mathbb{M}}
\newcommand{\bee}{\begin{equation}}
\newcommand{\ee}{\end{equation}}

\begin{document}

\cen{{\bf FRACALMO PRE-PRINT: \ http://www.fracalmo.org}}
% \vsh
\cen{{\bf Entropy (MDPI), Vol. 22  (2020) 1359/1--29; DOI:10.3390/e22121359}}
% \vsh
\cen{{\bf Special Issue: Fractional Calculus and the Future of Science}  }
%% \vsh
\cen{{\bf Guest Editor: Bruce West}}
\vsh
\hrule
% \end{center}
%%%%%%%%%%%%%%%%%%%%%%%%%%%%%%%%%%%%%%%%%%%%%%%%%%%%%%%%%%%%%%%%%%%%%%%%
\vskip 0.50truecm
\font\title=cmbx12 scaled\magstep1
\font\bfs=cmbx12 scaled\magstep1
\font\little=cmr10
\begin{center}
{\title Why the Mittag-Leffler function can be considered}
 \\ [0.25truecm]
 {\title the Queen function of the Fractional Calculus?}
 \\ [0.25truecm]
Francesco MAINARDI$^{(\star)}$ 
\\[0.25truecm]

$\null^{(\star)}$ {\little Department of Physics and Astronomy, University of Bologna, and INFN}
\\ 
{\little Via Irnerio 46, I-40126 Bologna, Italy}
\\{\little Corresponding Author.   E-mail: francesco.mainardi@bo.infn.it}
\\ %% \vskip 0.truecm
{Revised Version (December 2020) of the paper published in Entropy.}
%% (MDPI), 
%% with the addition of Appendix B: Some notes on the numerical computation of the %% functions of the Mittag-Leffler type.}

\end{center}
%%%%%%%%%%%%5
\begin{abstract}
%% Insert your abstract here.
\noindent 
In this survey we stress the importance of the higher transcendental 
Mittag-Leffler  function  in the framework of the Fractional Calculus.
We first start with the analytical properties of the classical Mittag-Leffler  
function as derived from being  the solution of the simplest fractional differential equation
 governing relaxation processes.
Through the Sections of  the text we plan  to address
 the reader in this pathway towards the main applications of the Mittag-Leffler  function that has induced us in the past to define it as the {\it Queen Function of the Fractional Calculus}.
 These applications concern some noteworthy  stochastic processes and the time fractional diffusion-wave equation.
 We expect that in the next future this function will gain more credit in the science of complex systems.
 In Appendix A we sketch some historical aspects related to the author's acquaintance
 with this function. Finally, with respect to the published version in Entropy, 
  we add Appendix B where we briefly refer to the numerical methods nowadays available to compute the functions of the Mittag-Leffler type.   
\end{abstract}
%
% \maketitle

\vsp
{\it 2010 Mathematics Subject Classification (MSC)}:
26A33, 33E12,  44A10.
%% 26A33,  %%%%  (main);    Fractional derivatives and integrals
%% 33E12, %% Mittag-Leffler type functions
%% 44A10,  %% Laplace Transforms
\vsp
{\it Key Words and Phrases}: 
Fractional Calculus; Mittag-Lefflller Functions; Wright Functions; 
 Complete Monotonicity;
 Fractional Differential Equations;
Fractional Relaxation; Diffusion-Wave Equation; Laplace and Fourier Transforms;
Fractional Poisson Process; Continuous Time Random Walks;  Complex Systems

\section{Introduction}

%\vsp  % please note that we have indent here. 
For  a few decades, the special transcendental function  known as the Mittag-Leffler function has attracted
the increasing  attention of  researchers  because of its key role in treating problems related to integral and differential equations of fractional order.

%\vsp
Since its introduction in 1903--1905 by the Swedish mathematician Mittag-Leffler at the beginning of the last century up to the 1990s,  this function was seldom considered by  mathematicians and applied scientists.

%\vsp
Before the 1990s, from a mathematical point of view, we recall the
1930  paper by Hille and Tamarkin \cite{Hille-Tamarkin 1930}
on the solutions of the Abel integral equation of the second kind,
and the  books by Davis \cite{Davis BOOK36},
Sansone \&  Gerretsen \cite{Sansone-Gerretsen BOOK60},
Dzherbashyan \cite{Dzherbashyan BOOK66} (unfortunately in Russian),
and finally  Samko et al. \cite{SKM BOOK93}.
Particular mention would be for
 the 1955 Handbook of  High Transcendental Functions of the Bateman project 
 \cite{Erdelyi BATEMAN},
 where this function  was  treated in Vol. 3, in the chapter devoted to miscellaneous functions.
For former applications we recall
an interesting  note by Davis \cite{Davis BOOK36} reporting previous research   by Dr. Kenneth S. Cole  in connection with nerve conduction,
and the  papers    by  Cole \& Cole \cite{Cole-Cole 1942},
Gross \cite{Gross JAP47}  and
Caputo \& Mainardi
\cite{Caputo-Mainardi PAGEOPH1971,Caputo-Mainardi RNC1971},
where the Mittag-Leffler function was adopted to represent
the responses in dielectric and   viscoelastic media. More information 
are found in the Appendix of this survey.

%\vsp
In the 1960's  the Mittag-Leffler function started to exit from the realm of miscellaneous functions because it was considered as a special case of the general  class of  Fox $H$ functions, that
can exhibit an arbitrary number of  parameters in their integral Mellin-Barnes representation. %see e.g.the books by
%Kiryakova \cite{Kiryakova_BOOK1994},
% Kilbas and Saigo  \cite{Kilbas-Saigo_BOOK2004},
%Marichev  \cite{Marichev_BOOK1983},
% Mathai \& Saxena \cite{Mathai-Saxena_BOOK1978},
% Mathai et al. \cite{Mathai-Saxena-Haubold_BOOK2010},
% Srivastava et.al. \cite{Srivastava-Gupta-Goyal_BOOK1982}.
%\vsp
However, in our opinion, this classification in  a too general framework has,  to some extent, obscured the relevance and the applicability of this function in applied sciences.
In fact,  most mathematical models are based  on a small  number of parameters, say  1 or 2 or 3,
so that a general theory may be confusing whereas the adoption of a generalized Mittag-Leffler function
with 2 or 3 indices may be sufficient.
% see e.g.
%Beghin \& Orsingher \cite{Beghin-Orsingher_EPJ2010},
%Capelas et al. \cite {Capelas-et-al_EPJ-ST2011},
%  Sandev et al. \cite{Sandev_FLE_FCAA2012},
% Tomovski et al. \cite{Tomovski-et-al_ITSF2010}.
%% \vsp
%Multi-index Mittag-Leffler functions have been introduced as well,
%see e.g.
%Kilbas et al. \cite{Kilbas-et-al_FCAA2013},
%Kilbas \& Saigo \cite{Kilbas-Saigo_1995},
%Kiryakova \cite{Kiryakova_CMA2010},
% Kiryakova \&  Luchko \cite{Kiryakova-Luchko_2010},
% but their extensive use has  not yet been  pointed out in applied sciences
%  up to now.

%\vsp
Nowadays it is well recognized that the Mittag-Leffler function  
 plays  a fundamental role  in  Fractional Calculus
%% in the complex and real domains and also in
%% the fields of Integral trasnsform (specially Laplace, Fourier, Mellin)
even if with a single parameter (as originally introduced by Mittag-Leffler)  just to be worth of being  referred to as  the {\it Queen Function of Fractional Calculus}, see Mainardi \&  Gorenflo~\cite{Mainardi-Gorenflo FCAA07}.
We find some information on the Mittag-Leffler functions in any treatise on Fractional Calculus
but for more details we refer the reader to the surveys of Haubold, Mathai and Saxena
\cite{Haubold-Mathai-Saxena 2011} and by Van Mieghem
\cite {Van Mieghem arXiv2020} and to the treatise by 
Gorenflo et al. \cite{GKMR BOOK20},
just devoted to Mittag-Leffler functions, related topics and applications. 

%\vsp 
The plan of this survey is the following. We start to give in Section 2 the main definitions and properties of the Mittag-Leffler function in one parameter 
  with related Laplace transforms.
Then in Section 3 we describe its use in the simplest
 fractional relaxation equation
 pointing out its  compete monotonicity.
 The asymptotic properties are briefly discussed in Section 4.
 In section 5 we briefly discuss the so called generalized Mittag-Leffler function, that is 
 the 2-parameter Mittag-Leffler function. Of course further generalization to 3 and more parameter will be referred to specialized papers and books.  
 Then in the following sections we discuss the application of the Mittag-Leffler function  in some noteworthy stochastic processes.
 We start in Section 6 with  the fractional Poisson process, and then in Section 7 with its application of the thinning of  renewal processes.
 The main application are dealt in Section 8 where we discuss the continuous time random walks (CTRW) and then in Section 9  we point out the asymptotic universality.
  In Section 10 we discuss the  time fractional diffusion-wave processes pointing out the role of the Mittag-Leffler functions in two parameters and their connection with the basic Wright functions. 
  In Section 11 we  draw our conclusive remarks.

The survey is completed with two additional appendices.
%%% 
 In Appendix A we find worthwhile to report the acquaintance of the author 
 with the Mittag-Leffler functions started in the late 1960s and 
 continued up to nowadays.
 Finally, we add Appendix B where we briefly refer to the numerical methods nowadays available to compute the functions of the Mittag-Leffler type.
 
% \vsp 
 We recall that Sections 3--10 are taken from  several papers by the author,
 published alone and with colleagues and former students.    
 Furthermore we have not considered other applications of the Mittag-Leffler functions including, for example,  anomalous diffusion theory in terms of fractional and generalized Langevin equations. On this respect we refer the readers to the articles of the author,
 see References \cite{Mainardi-Pironi EM96}, \cite{Mainardi-Mura-Tampieri MPS09}, and to the 
 recent book by Sandev and Tomovski \cite{Sandev-Tomovski BOOK19} and references therein.
 For many  items related to  the Mittag-Leffler functions we refer
 again to the treatise by Gorenflo et al. \cite{GKMR BOOK20}.

%On this respect  we  just point out  some
% reviews and treatises on Fractional Calculus  (in order of publication time)
% that have more  attracted our attention:
%  Gorenflo \& Mainardi \cite{Gorenflo-Mainardi CISM97},
%  Podlubny \cite{Podlubny BOOK99},
%  Hilfer \cite{Hilfer BOOK00},
%  Kilbas et al. \cite{Kilbas-Srivastava-Trujillo BOOK206},
%  Magin \cite{Magin_BOOK06},
%  Mathai \& Haubold \cite{Mathai-Haubold BOOK08},
%  Mainardi \cite{Mainardi BOOK10},
%  Diethelm \cite{Diethelm BOOK10},
%  Tarasov \cite{Tarasov BOOK11},
%   Klafter et al. \cite{KLM BOOK12},
%   Baleanu et al. \cite{BDST BOOK12},
%   Uchaikin \cite{Uchaikin BOOK13},
%  Capelas \cite{Capelas BOOK19}
%   along with the treatise by Gorenflo et al. \cite{GKMR BOOK20} 
%   devoted to Mittag-Leffler functions, related topics and applications. 
 %% \vsp
  % \newpage
\section{The Mittag-Leffler functions: definitions and Laplace transforms. }
%%\noindent
The Mittag-Leffler function is defined by the following power series,
convergent in the whole complex plane,
%%$$
\bee
E_\alpha (z) := \sum_{n=0}^\infty \frac{z^n}{\Gamma (\alpha n+1)}
\,,\q \alpha > 0\,, \q z\in \CC\,. 
\label{2.1} %%\eqno(2.1)$$
\ee
%% so $E_\alpha (z) $ is an {\it entire function}.
We recognize that it is an entire function of 
order $1/\alpha$  providing  a simple generalization of the
  exponential function $\exp (z)$  to which it reduces for $\alpha=1$.
For detailed information on the Mittag-Leffler-type functions
and their Laplace transforms the reader  may  consult \eg
\cite{Erdelyi BATEMAN,Gorenflo-Mainardi CISM97,Podlubny BOOK99} and
 the recent treatise by Gorenflo et al. \cite{GKMR BOOK20}.

 %%
 %\vsp
 %REMINDING  (\ref{2.1}) WE GET
 %\vsp
  We also note that for the convergence of the power series in (\ref{2.1})
 the parameter $\alpha$ may be complex provided that $\Re (\alpha) >0$.
%% as pointed out in \cite{Mittag-Leffler 04}.
%\vsp
The most interesting properties of the Mittag-Leffler function
are associated with its asymptotic expansions
as $z \to \infty$ in various sectors of the complex plane.
% \vsp
%For  detailed asymptotic analysis,
%which includes the smooth transition across the Stokes lines,
% the interested reader is referred to Wong and Zhao \cite{Wong-Zhao_CA2002} and to 
% Paris \cite{Paris ML}.
\vsp
In this paper  we mainly consider  the Mittag-Leffler function of order $\alpha \in (0,1)$ on the negative real semi-axis where is known to be completely monotone (CM) due a classical  result by  Pollard \cite{Pollard BAMS48}, see also
Feller \cite{Feller BOOK71}.
%%  who proved a conjecture formerly pointed out by  Feller.
 % (based on the probability theory).

%\vsp
Let us recall that a function $ \phi(t)$  with $t\in\RR^+$ is    called a  completely
monotone (CM) function if it is non-negative, of class $C^{\infty}$,  and
$(-1)^n \phi^{(n)}(t)\ge 0$ for all $n \in \NN$.
Then a function $ \psi(t)$  with $t\in\RR^+$ is    called a Bernstein function
if it is non-negative, of class $C^{\infty}$, with a CM first derivative.
These functions play fundamental roles in linear hereditary mechanics to represent relaxation and creep processes, see, for example, Mainardi \cite{Mainardi BOOK10}.
For mathematical details we refer the interested reader  to the survey paper
by Miller and  Samko \cite{Miller-Samko ITSF01}
and to the treatise by Schilling et al. \cite{SSV BOOK12}.

%\vsp
In particular we are interested in the function
\bee   %%5$$ 
e_\alpha(t) := E_\alpha(-t^\alpha)
 =
\sum_{n=0}^\infty (-1)^n \frac{t^{\alpha n}}{\Gamma (\alpha n+1)}
\,, \quad t>0\,, \quad 0<\alpha \le 1\,,
\label{2.2} 
\ee %%\eqno (2.2)$$
whose Laplace transform pair reads
\bee %%$$
       e_\alpha(t) \,\div \, \frac{s^{\alpha-1}}{s^\alpha +1}\,, \quad \alpha >0 \,. 
  \label{2.11} %%     \eqno(2.11)$$
  \ee
  Here we have used the notation $\div$ to denote the juxtaposition of a function of time  
  $f(t)$ with its Laplace transform
  $$ \widetilde {f}(s) = \int _0^\infty \e^{-st}\, f(t)\, dt\,. $$
The pair (\ref{2.11})    
  can be proved   by transforming term by term the power series representation  
  of $ e_\alpha(t)$
 in the R.H.S of  (\ref{2.2}).
Similarly we can prove the following Laplace transform pair for its time derivative 
\bee %%$$
       e^\prime_\alpha(t)=
\frac{d}{dt} E_\alpha(-t^\alpha)       
        \,\div \, - \frac{1}{s^\alpha +1}\,, \quad \alpha >0 \,. 
  \label{2.01} %%     \eqno(2.11)$$
  \ee
 For this Laplace transform pair  we can simply apply the usual rule for the Laplace transform 
  for the first derivative of a function, that reads  
  $$ \frac{d}{dt} f(t) \, \div \, s\, \widetilde{f}(s) - f(0^+)\,. $$

\vsp
\section{The Mittag-Leffler function in fractional relaxation processes}

For readers'  convenience let us briefly outline the topic
concerning the   generalization  via fractional calculus
 %% of to the initial value problem for
 of the first-order differential equation governing the
 phenomenon of (exponential) relaxation.
Recalling   (in non-dimensional units)
 the initial value problem %% standard equation
 \bee %% $$
 \frac{du}{dt} = -u(t) \,, \quad t\ge 0\,, \quad \hbox{with}\quad u(0^+)= 1\,
 \label{2.3} %%\eqno(2.3)$$
 \ee
 whose solution is
 %% The well-known solution of the above initial value problem
\bee %% $$
u(t) = \exp (-t)\,,
\label{2.4} %%\eqno (2.4)$$
\ee
  the following  two alternatives %%  with respect to the R-L and C fractional derivatives
  with $\alpha \in (0,1)$ are offered in the literature:
  \bee %% $$
 (a)\quad \frac{du}{dt } = - D_t^{1-\alpha}\,u(t) \,,\quad  t\ge 0\,, \quad \hbox{with}\quad u(0^+)= 1\,,
  \label{2.5a} %%\eqno(2.5a)$$
  \ee
 \bee%% $$ 
 (b)\quad  _*D_t^\alpha \,u(t) = -u(t)\,, \quad t\ge 0\,, 
 \quad \hbox{with}\quad u(0^+)= 1\,,
 \label{2.5b}  %%\eqno(2.5b)$$
 \ee
where
 $D_t^{1-\alpha}$ and $\,  _*D_t^\alpha$ denote the fractional derivative of order $1-\alpha$ in
 the Riemann-Liouville sense and the fractional derivative of order $\alpha$ in the Caputo sense, respectively.
 \vsp
 For a generic order $\mu\in (0,1)$ and for a sufficiently well-behaved function $f(t)$ with $t\in \RR^+$
 the above derivatives are defined as follows, see for example,
 Gorenflo and Mainardi \cite{Gorenflo-Mainardi CISM97},
 Podlubny \cite{Podlubny BOOK99},
 \bee  %% $$ 
 (a) \quad D_t^\mu  \,f(t) =
  {\ds  \rec{\Gamma(1-\mu )}}\,
{\ds \frac{d}{dt}}\left[
  \int_0^t
    \! \frac{f(\tau)}{ (t-\tau )^{\mu }}\,d\tau\right]\,, 
   \label{2.6a} %% \eqno(2.6a)  $$
   \ee
   \bee %% $$  
 (b) \quad  _*D_t^\mu  \,f(t) =
{\ds \rec{\Gamma(1-\mu )}}\int_0^t
    \! \frac{f^\prime (\tau)}{ (t-\tau )^{\mu }}\,d\tau\,. 
    \label{2.6b} %%\eqno(2.6b)  $$
    \ee
  Between the two derivatives we have the  relationship
  \bee %% $$
	 {\ds _*D_t^\mu  \,f(t)} =
  {\ds D_t^\mu \,f(t) -  f(0^+)\, \frac{t^{-\mu}} {\Gamma(1-\mu )}}  =
 D_t^\mu  \,\left[ f(t) - f(0^+) \right] \,.
 \label{2.7} %%\eqno(2.7}  $$
 \ee
Both derivatives in the limit $\mu \to 1^-$ reduce to the standard first derivative but for
$\mu \to 0^+$ we have
\bee %% $$ 
D_t^\mu f(t) \to f(t)\,, \quad  _*D_t^\mu f(t) = f(t) - f(0^+)\,, \quad \mu \to 0^+\,.
\label{2.8}%%\eqno(2.8)$$
\ee
 %\vsp
 In analogy to the standard problem (\ref{2.3}), we  solve the problems 
 (\ref{2.5a}) and (\ref{2.5b})
  with the Laplace transform technique,  using  the rules pertinent to the corresponding fractional derivatives, that we recall hereafter for a generic order $\mu \in (0,1)$,
   \bee  %% $$ 
 (a) \quad D_t^\mu  \,f(t) \, \div s^\mu \widetilde{f}(s) - g(0^+)\,,
 \quad g(0^+)= \frac{1}{\Gamma(1-\mu)}\, \lim_{t\to 0^+} \int_0^t 
 (t-\tau)^{-\mu}\, f(\tau)\, d \tau\,.
 \label{LT.DRL}
 \ee
   \bee  %% $$ 
 (b) \quad _*D_t^\mu  \,f(t) \, \div s^\mu \widetilde{f}(s) - f(0^+)
 \,.
 \label{LT.DC}
 \ee
 We note that it is generally more cumbersome to use the Laplace transform pair for the Rieaman Liouville derivative (\ref{LT.DRL}) that for the Capute derivative (\ref{LT.DC}). Indeed 
 the rule (\ref{LT.DRL} requires the initial value of the fractional integral of $f(t)$
 whereas the rule (\ref{LT.DC})  simply requires the initial value of $f(t)$.
 For this property the Caputo derivative is mostly used in physical problems where finite  
 initial values are given. 

% \vsp 
 Then we recognize that  the problems (a) and (b) are equivalent since  the Laplace transform of the solution
 in both cases comes out as
 \bee %% $$ 
  \widetilde u(s) = \frac{s^{\alpha-1}}{s^\alpha +1}\,, 
  \label{2.9} %%\eqno(2.9) $$
  \ee
that yields, in virtue of the Laplace transform pair (\ref{2.11}) 
  \bee %% $$ 
u(t) = e_\alpha(t) := E_\alpha(-t^\alpha)\,.
\label{2.10} %%\eqno(2.10)$$
\ee
We thus recognize that the Mittag-Leffler function  provides the solution to the fractional relaxation equation, as outlined, for example, by
Gorenflo and Mainardi \cite{Gorenflo-Mainardi CISM97},
Mainardi and  Gorenflo \cite{Mainardi-Gorenflo FCAA07},
and Mainardi \cite{Mainardi BOOK10}.

%\vsp
 Furthermore,    by
 anti-transforming the R.H.S  of (\ref{2.11})  by using the complex Bromwich formula, and taking into account  for $0<\alpha<1$ the contribution from  branch cut on the negative real semi-axis (the denominator $s^\alpha +1$  does nowhere vanish in the cut plane
  $-\pi \le \hbox{arg} \, s \le \pi $),  we get,
 see the survey by  Gorenflo and Mainardi \cite{Gorenflo-Mainardi CISM97},
\bee  %% $$
 e_\alpha(t) = \int_0^\infty \!\!  \e^{-rt} K_\alpha(r)\, dr\,,
 \label{2.12} %% \eqno(2.12) $$
 \ee
where
\bee %$$ 
     K_\alpha(r) =
 \mp \,\rec{\pi}\,    {\rm Im}\,
  \left\{ \left. \frac{s^{\alpha -1}} {s^\alpha +1}\right\vert_{{\ds s=r\, \e^{\pm i\pi}}} \right\}
 = \rec{\pi}\,
   \frac{ r^{\alpha -1}\, \sin \,(\alpha \pi)}
    {r^{2\alpha} + 2\, r^{\alpha} \, \cos \, (\alpha \pi) +1} \ge 0\,.
\label{2.13} %%      \eqno(2.13)$$
\ee
%%%%%
We note that this formula was obtained as a simple exercise of complex analysis
without be aware of  the Titchmarsh  formula for inversion of Laplace transforms \cite{Titchmarsh BOOK37},
revised  by  Gross and Levi \cite{Gross-Levi MATHNOTES46} and by Gross 
\cite{Gross PhilMag50}.
This formula is rarely outlined  in books on Laplace transforms so we refer the reader
for example  to    Apelblat's book \cite{Apelblat BOOK11} for its presence.
\vsp
    Since   $K_\alpha(r)$ is non-negative for all $r$ in the integral, the above formula proves
    that $e_\alpha(t)$ is a
     CM function in view of the Bernstein theorem. 
     This theorem provides a necessary and sufficient condition for a CM function as a real Laplace transform of a non-negative measure.

%\vsp
However, the CM property of $e_\alpha(t)$  can also be seen  as a consequence  of the result by  Pollard~\cite{Pollard BAMS48}
because the transformation $x=t^\alpha$ is a Bernstein function for $\alpha \in (0,1)$.
In fact it is known that a CM function can be obtained by composing a CM with a Bernstein function based on the following theorem:
 %% \vsp
{\it Let $\phi(t)$ be a CM function and let $\psi(t)$ be a Bernstein function, then $\phi[\psi(t)]$  is a CM function.}
 \vsp
  As a matter of fact,  $K_\alpha(r)$ provides an interesting  spectral representation of  $e_\alpha(t)$ in frequencies. With the change of variable $\tau=1/r$ we get the corresponding spectral representation in relaxation times, namely
\bee %$$
 e_\alpha(t) = \int_0^\infty \e^{-t/\tau} H_\alpha(\tau)\, d\tau\,, \;
H_\alpha(\tau) = \tau^{-2} \, K_\alpha(1/\tau)\,, 
\label{2.14} %%\eqno(2.14) $$
\ee
that can be interpreted as a continuous distributions of elementary (i.e., exponential) relaxation processes.
As a consequence    we get the identity between  the two spectral distributions, that is
\bee %% $$ 
K_\alpha(r)  = H_\alpha(\tau)
= \rec{\pi}\,
   \frac{ \tau^{\alpha -1}\, \sin \,(\alpha \pi)}
    {\tau^{2\alpha} + 2\, \tau^{\alpha} \, \cos \, (\alpha \pi) +1}
\,, 
\label{2.15} %%\eqno(2.15)$$
\ee
a surprising fact pointed out in Linear Viscoelasticity
by the author in his  book \cite{Mainardi BOOK10}.
This kind of universal/scaling  property seems a peculiar one for our Mittag-Leffler function $e_\alpha(t)$.

\begin{figure}  %%[H]
\centering
\includegraphics[width=7.5cm]{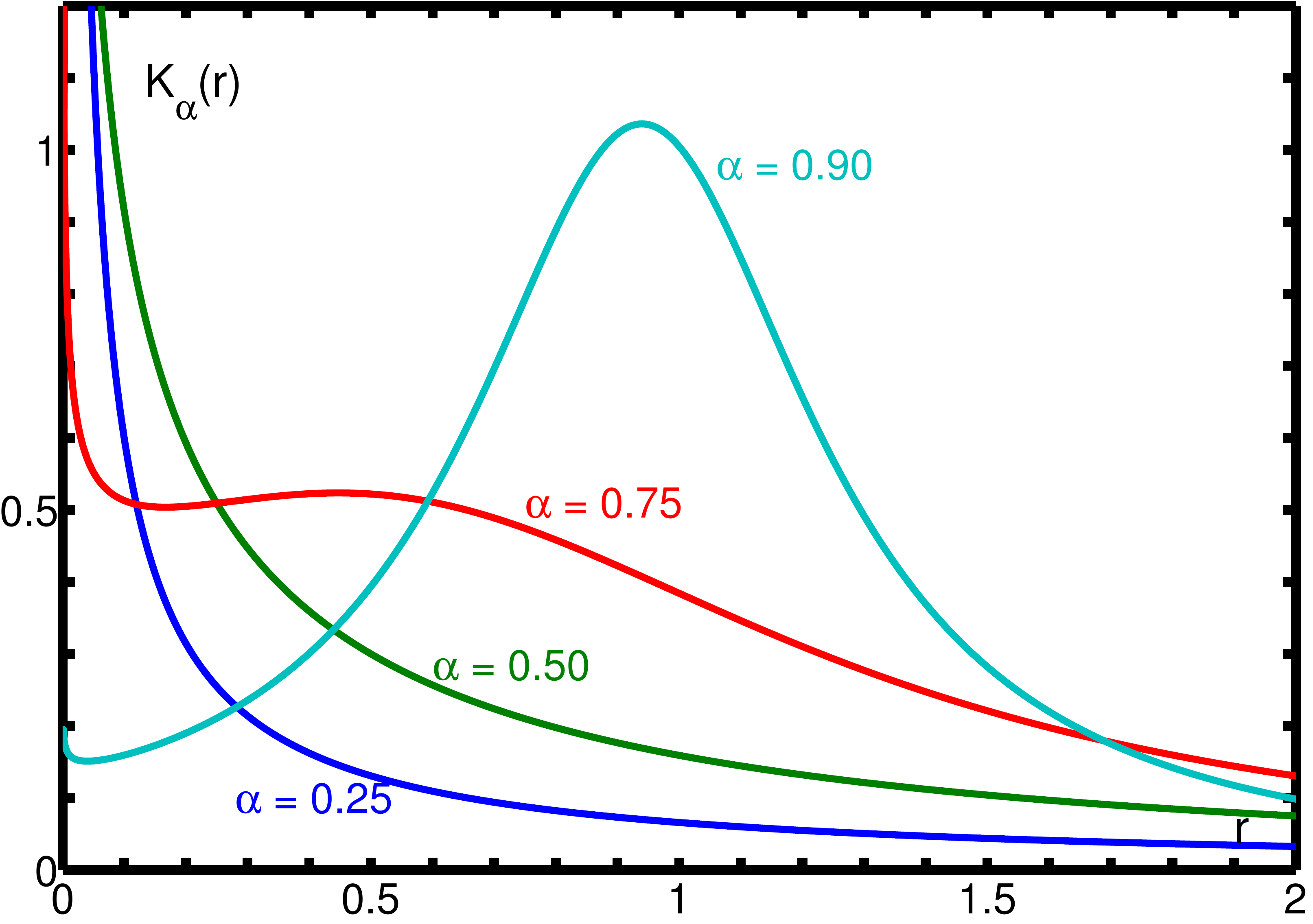}
\caption{The spectral function  $K_\alpha(r) $ for $\alpha=0.25, 0.50, 0.75, 0.90$ in the frequency range  $0\le r \le 2$.}
\label{fig1}
\end{figure}

%\vsp
In Figure~\ref{fig1}, we show $K_\alpha(r)$ for some values of the parameter $\alpha$.
 Of course for $\alpha=1$ the Mittag-Leffler function reduces to the exponential function 
 $\exp(-t)$
 and the  corresponding spectral distribution is the    Dirac delta generalized function centered at $r=1$, namely $\delta(r-1)$.
\vsp

% \begin{center}
%\includegraphics[width=7.5cm]{MDC_fig2a}
%\end{center}
%%\caption
%{{\bf Fig.1} The spectral function  $K_\alpha(r) $ for $\alpha=0.25, 0.50, 0.75, 0.90$ in the frequency range  $0\le r \le 2$.}
 %%%%%%%
 %\newpage
 %%%%

\begin{figure} %% [H]
\centering
\includegraphics[width=7.5cm]{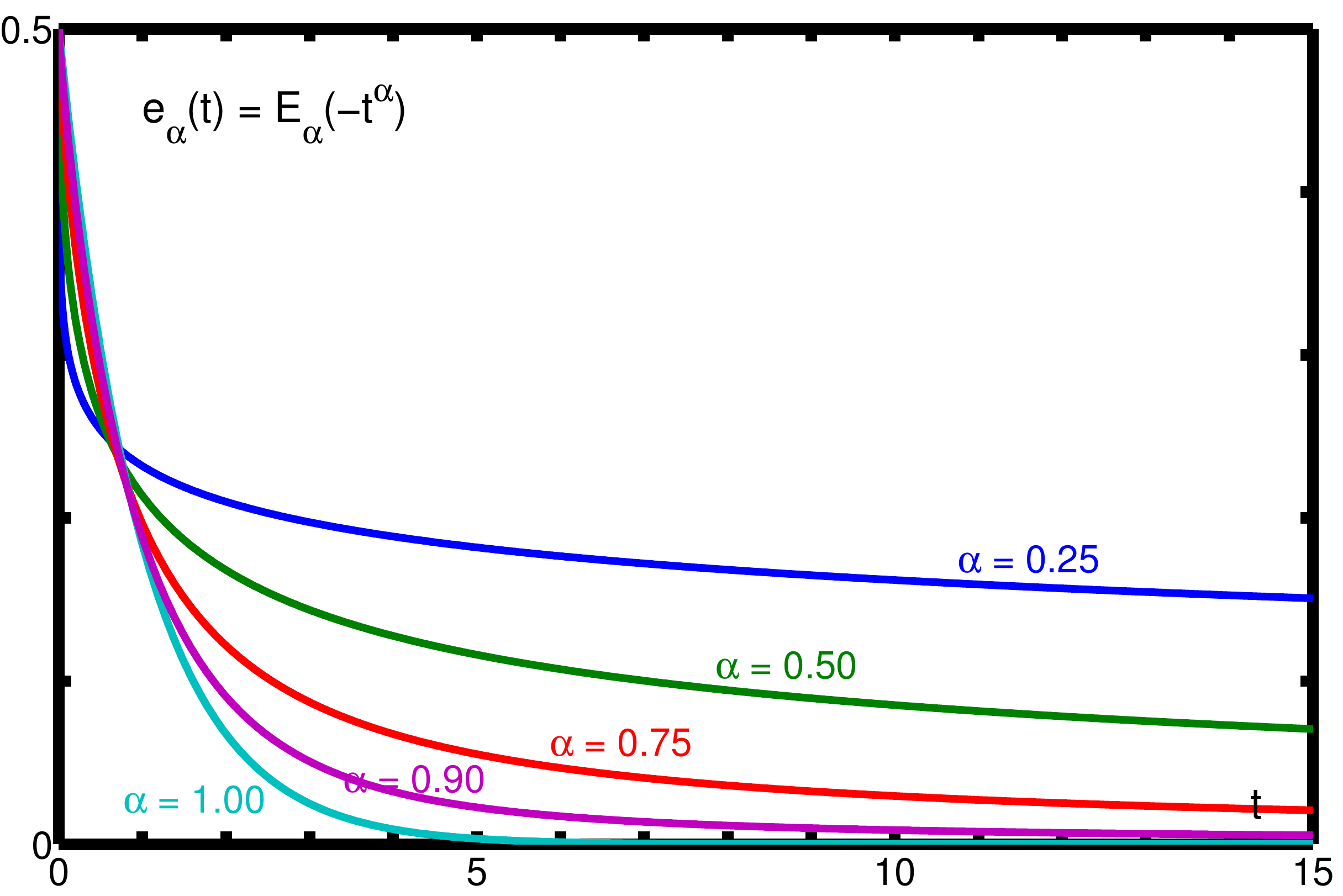}
\caption{The Mittag-Leffler function $e_\alpha(t)$  for $\alpha=0.25, 0.50, 0.75, 0.90, 1.$ in the time range $0\le t\le 15$.}
\label{fig2}
\end{figure}

%\vsp
In Figure~\ref{fig2},  we show some plots of $e_\alpha(t)$  for some values of the parameter $\alpha$.
It is worth to note the different rates of decay of $e_\alpha(t) $   for small and large times.
In fact the decay is very fast as $t\to 0^+$ and very slow as $t \to +\infty$.
\vsp

%\begin{center}.
%\includegraphics[width=7.5cm]{MDC_fig3a_090}
%\end{center}
%%\caption
%{{\bf Fig.2} The Mittag-Leffler function $e_\alpha(t)$  for $\alpha=0.25, 0.50, 0.75, 0.90, 1.$ in the time range $0\le t\le 15$. }

%\vsp
The Mittag-Leffler function turns out the basic function in relaxation processes of physical interested as occurring in viscoelastic and dielectric materials.
We refer the readers for viscoelasticity, that is, to the contribution of the author including 
References \cite{Mainardi BOOK10},\cite{Mainardi HFCA19} and 
Reference \cite{Mainardi-Spada EPJ-ST11} 
whereas for dielectric 
materials to the survey by Garrappa et al. 
\cite{Garrappa-Mainardi-Maione FCAA16}.
For the pioneers who have pointed out  the role of the Mittagf-Leffler function 
in mechanical and dielectric relaxation processes we refer to the recent survey by Mainardi and Consiglio \cite{Mainardi-Consiglio WSEAS2020}  

\section{Asymptotic approximations to the Mittag-Lefler function}

We now report the two common asymptotic approximations of our Mittag-Lefflrer function.
 Indeed,  it is common to point out that     the function $e_\alpha(t)$  matches for $t\to 0^+$ with a stretched exponential   
  with an infinite negative derivative,   whereas as $t\to \infty$ 
  with a negative power law.
The short time  approximation is  derived from the convergent power series representation (\ref{2.2}). 
In  fact,
\bee %% $$
 e_\alpha(t) = 1 - \frac{t^\alpha}{\Gamma(1+\alpha)} + \dots
\sim  \exp{\ds \left[- \frac{t^\alpha}{\Gamma(1+\alpha )}\right]}\,, \quad  t\to 0\,.
\label{3.1} %% \eqno(3.1)$$
\ee
The long time approximation is derived from the asymptotic power series representation of $e_\alpha(t)$ that turns out to be, see 
\cite{Erdelyi BATEMAN} %% Erd\'elyi (1955),
\bee %%  $$ 
e_\alpha(t) \sim
 \sum_{n=1}^\infty  (-1)^{n-1} \,\frac{ t^{-\alpha n}}{\Gamma(1- \alpha n)}\,,
 \quad t\to \infty\,,
 \label{3.2} %%\eqno(3.2) $$
 \ee
 so that, at the first order,
\bee %% $$
e_\alpha(t) \sim  {\ds \frac{t^{-\alpha}}{\Gamma(1-\alpha )} }\,, \quad t\to \infty\,. 
\label{3.3} %%\eqno(3.3) $$
\ee
 As a consequence  the function $e_\alpha(t)$  interpolates
 for intermediate time $t$ between the stretched exponential
and the negative power law.
The stretched exponential models
  the  very fast decay  for small  time $t$, whereas the asymptotic   power law
  is due to the very slow decay for large  time $t$.
In fact, we have the two commonly stated  asymptotic  representations:
%% as $t\to 0$ and $t\to \infty$  commonly found in the literature:
\bee %%  $$
e_\alpha  (t) \sim
\left\{
\begin{array}{ll}
 %% 1- {\ds \frac{t^\alpha}{\Gamma(1+\alpha )} } \sim
 e_\alpha^0(t) := \exp{\ds \left[- \frac{t^\alpha}{\Gamma(1+\alpha )}\right]}\,,
  &  t\to 0\,;   \\ \\
e_\alpha^\infty(t) := {\ds \frac{t^{-\alpha}}{\Gamma(1-\alpha )} } =
{\ds \frac{\sin (\alpha  \pi)}{\pi}\,\frac{\Gamma(\alpha )}{t^\alpha }}\,,
 &  t\to \infty\,.
\end{array}
\right .
\label{3.4} %% \eqno(3.4) $$
\ee
   The stretched exponential   replaces  the rapidly decreasing expression
    $1- {t^\alpha}/{\Gamma(1+\alpha )} $  
    from (\ref{3.1}).
%%    asymptotically equivalent  to $e_\alpha(t)$ as $t\to 0$.
 Of course, {\it for sufficiently small and for sufficiently large values of $t$}
  we have the inequality
  \bee %%   $$ 
   e_\alpha^0(t) \le e_\alpha^\infty (t)\,,  \quad 0<\alpha<1\,. 
   \label{3.5} %%\eqno (3.5)$$
   \ee
   
\begin{figure} %% [H]
\centering
\includegraphics[width=6.35cm]{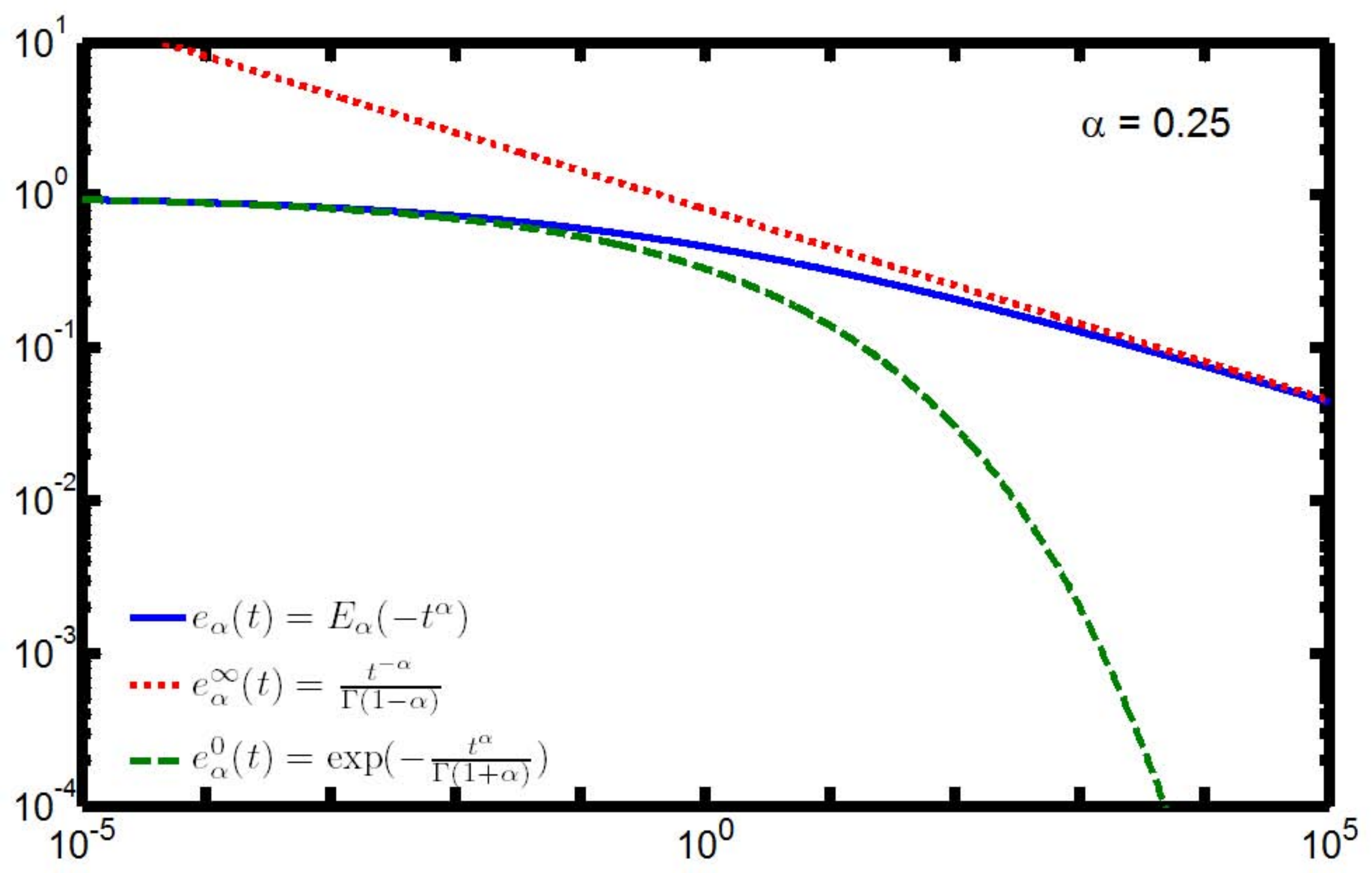}
\includegraphics[width=6.35cm]{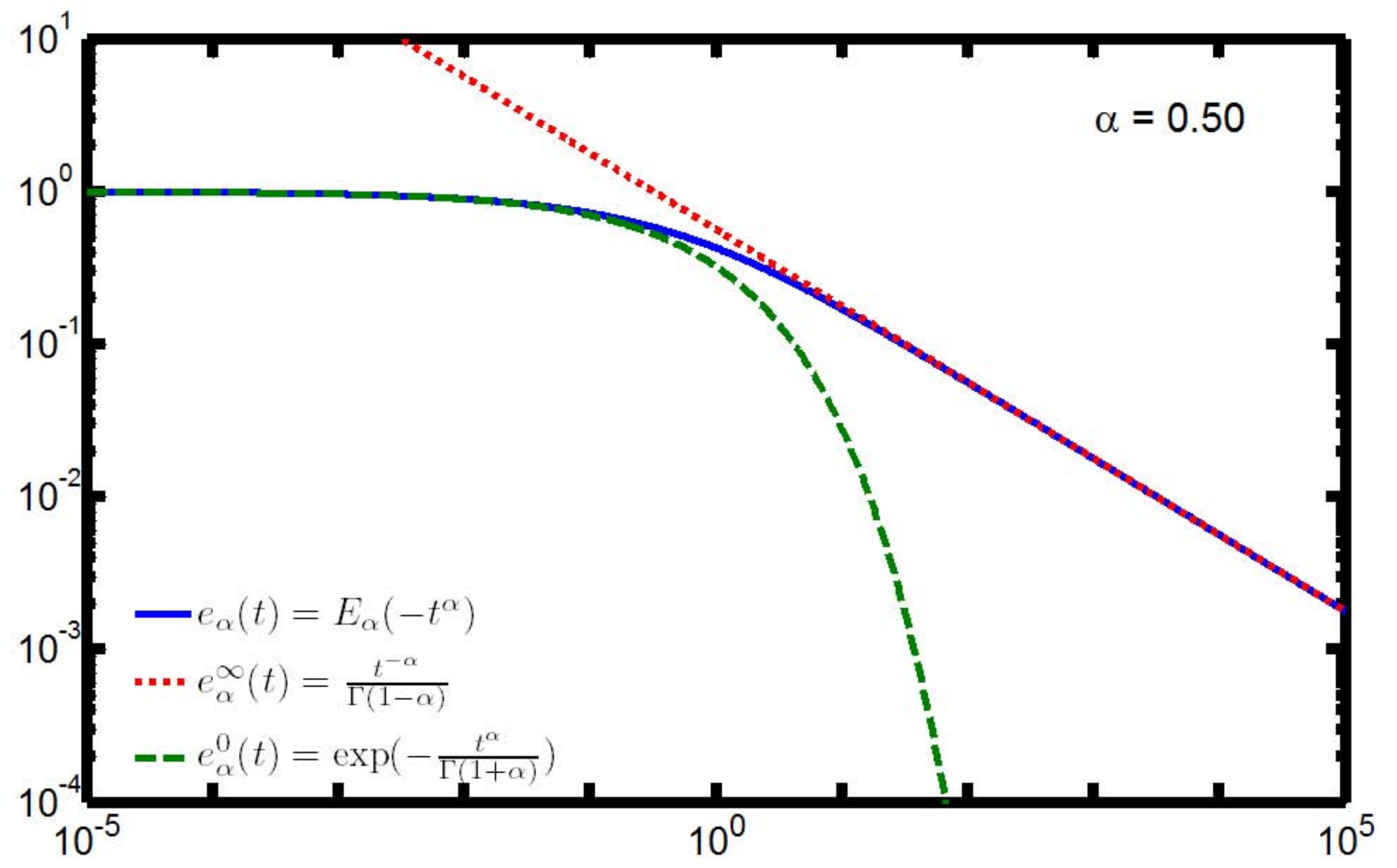}
\caption{Approximations $e^0_\alpha(t)$ (dashed line) and $ e^\infty_\alpha(t)$
(dotted line) to   $e_\alpha(t)$  
 in   $10^{-5} \le t \le 10^{+5}$ for $\alpha=0.25$ (LEFT)
 and for $\alpha=0/50$  (RIGHT).}
 \label{fig3}
\end{figure}

\begin{figure} %% [H]
\centering
\includegraphics[width=6.35cm]{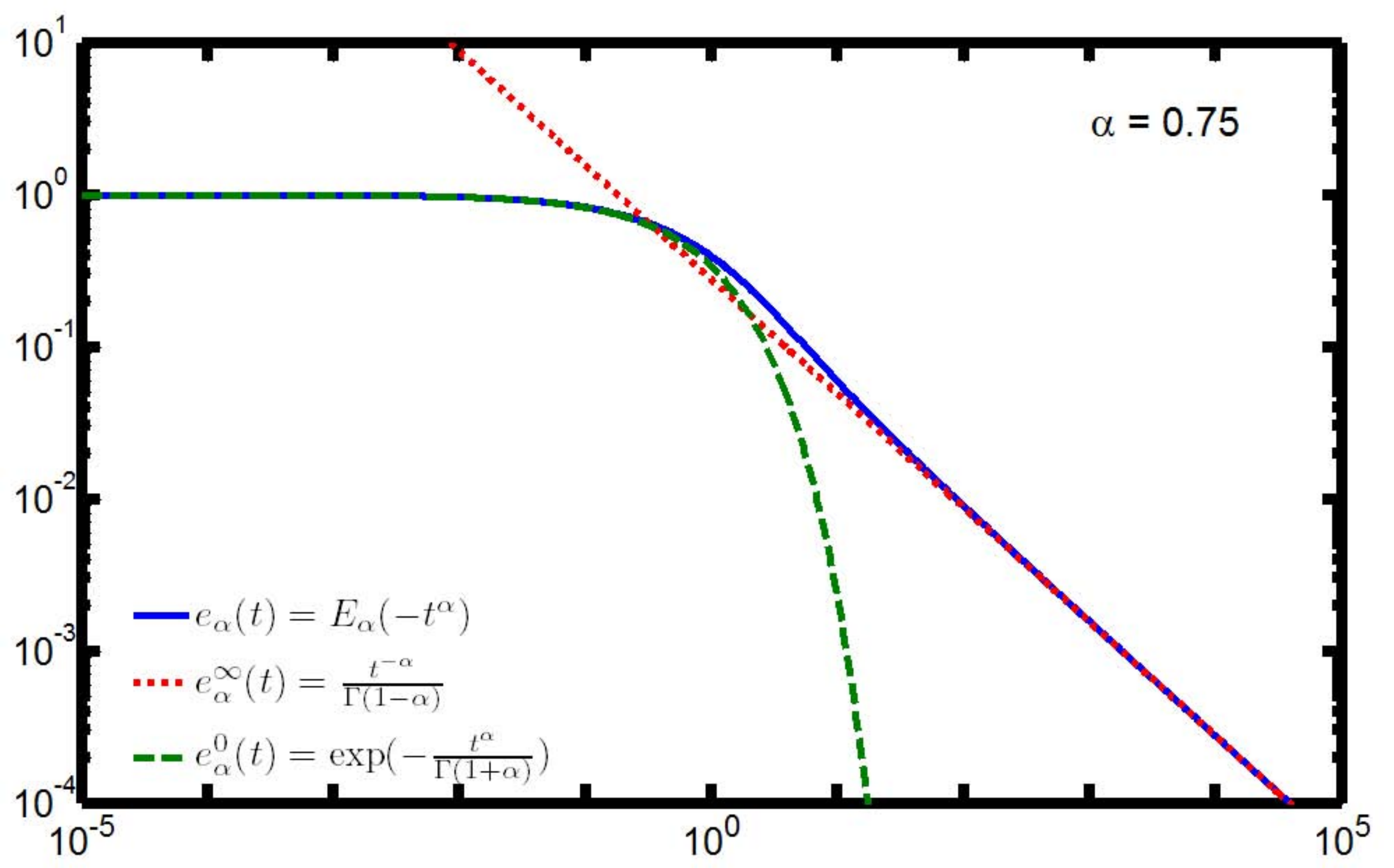}
\includegraphics[width=6.35cm]{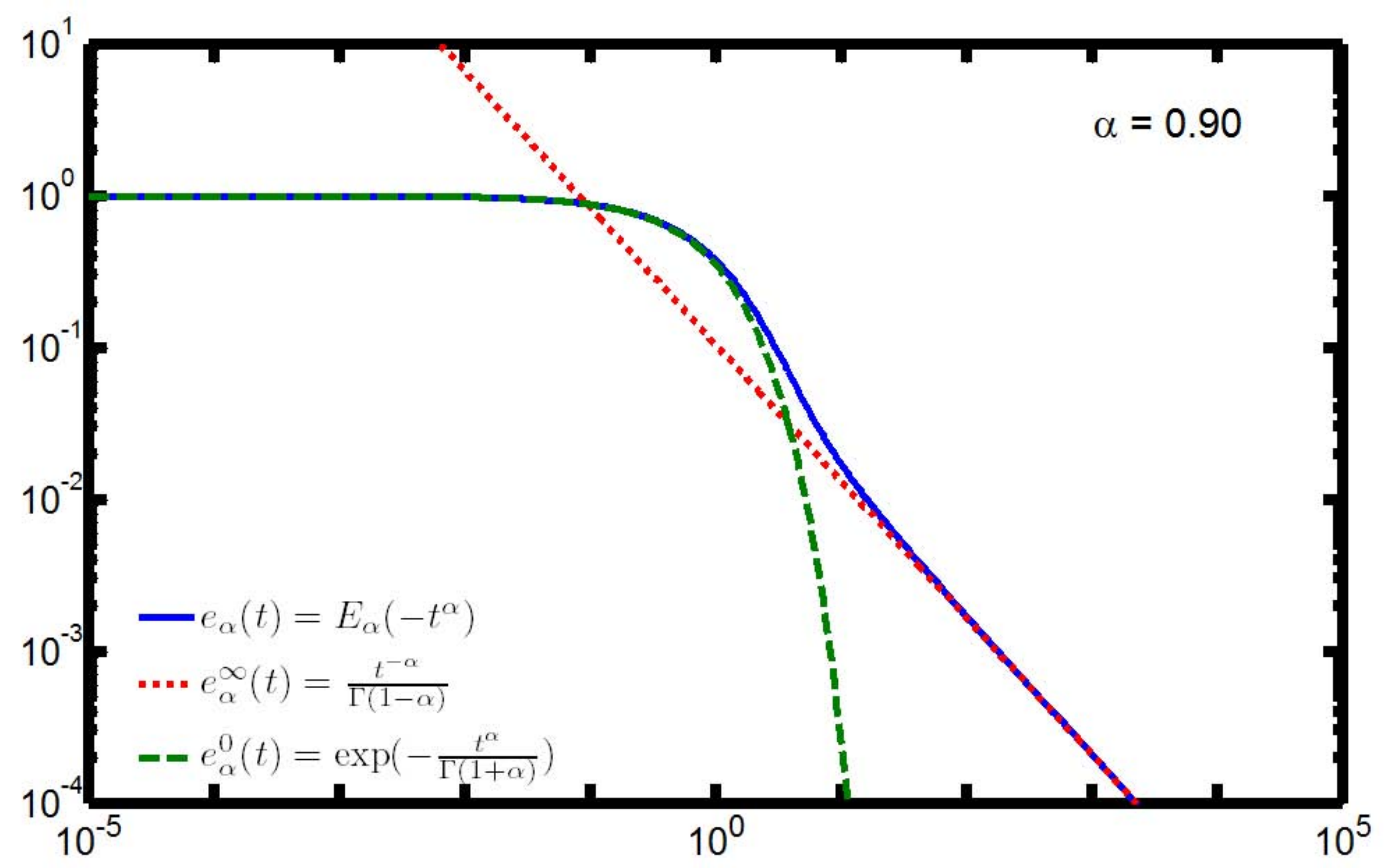}
\caption{Approximations $e^0_\alpha(t)$ (dashed line) and $ e^\infty_\alpha(t)$
(dotted line) to   $e_\alpha(t)$  (LEFT)
and  the corresponding relative errors (RIGHT)
 in   $10^{-5} \le t \le 10^{+5}$ for $\alpha=0/75$ (LEFT) and for
 $\alpha=0.90$ (RIGHT).}
 \label{fig4}
\end{figure}

%   \vsp
In  Figures~\ref{fig3} and \ref{fig4},  we compare for $\alpha = 0.25, 0.5, 0.75, 0.90$ in   logarithmic scales  the function $e_\alpha(t)$ (continuous line) and its asymptotic representations, the stretched exponential $e^0_\alpha(t)$ valid for $ t\to 0$ (dashed line)  and the power law $e^\infty_\alpha(t)$  valid for $t\to \infty$ (dotted line).
We have chosen the time range  $10^{-5} \le t \le 10^{+5}$.
 \vsp
We note from  these Figures
%% ~\ref{fig3} and \ref{fig4} 
that, whereas the plots of $e_\alpha^0(t)$ remain always under the corresponding ones of $e_\alpha(t)$, the plots of $e_\alpha^\infty(t)$ start above those of  $e_\alpha(t)$
but, at a certain point, an intersection may occur  so changing the sign of the relative errors.
The interested reader may consul the plots of the relative errors in the 2014 paper by 
the author %%Mainardi
\cite{Mainardi DCDS2014} from which, in particular, Figures 1--4 have been extracted.

%\begin{center}
%\includegraphics[width=6.35cm]{AG025.pdf}
%\includegraphics[width=6.35cm]{AG050.pdf}
%\end{center}
%%\caption
%{{\bf Fig.3} Approximations $e^0_\alpha(t)$ (dashed line) and $ e^\infty_\alpha(t)$
%(dotted line) to   $e_\alpha(t)$  
% in   $10^{-5} \le t \le 10^{+5}$ for $\alpha=0.25$ (LEFT)
% and for $\alpha=0/50$  (RIGHT).}
%%%%%%%%
%\smallskip
%%%%%%%%%%%%%
%%%%%%%%%%%%%
%\smallskip

%\begin{center}
%%\hskip -0.5truecm
%\includegraphics[width=6.35cm]{AG075.pdf}
%\includegraphics[width=6.35cm]{AG090.pdf}
%\end{center}
%% \caption
%{{\bf Fig.4} Approximations $e^0_\alpha(t)$ (dashed line) and $ e^\infty_\alpha(t)$
%(dotted line) to   $e_\alpha(t)$  (LEFT)
%and  the corresponding relative errors (RIGHT)
% in   $10^{-5} \le t \le 10^{+5}$ for $\alpha=0/75$ (LEFT) and for
% $\alpha=0.90$ (RIGHT).}
%%
%\newpage
% 

%%%%%%
\section{The generalized Mittag-Leffler function}
In this survey we devote our attention mainly to the classical Mittag-Leffler 
in one parameter $\alpha$  as introduced by Mittag-Leffler since 1903 defined by the power series  in (\ref{2.1}).
We have just learned from the instructive E-print by Van Mieghem 
\cite{Van Mieghem arXiv2020} that
the series (\ref{2.1})  was discussed by Hadamard in 1893, that is  
10 years earlier  than Mittag-Leffler himself.

\begin{figure}   %% [H]
\begin{center}
\includegraphics[width=.80\textwidth]{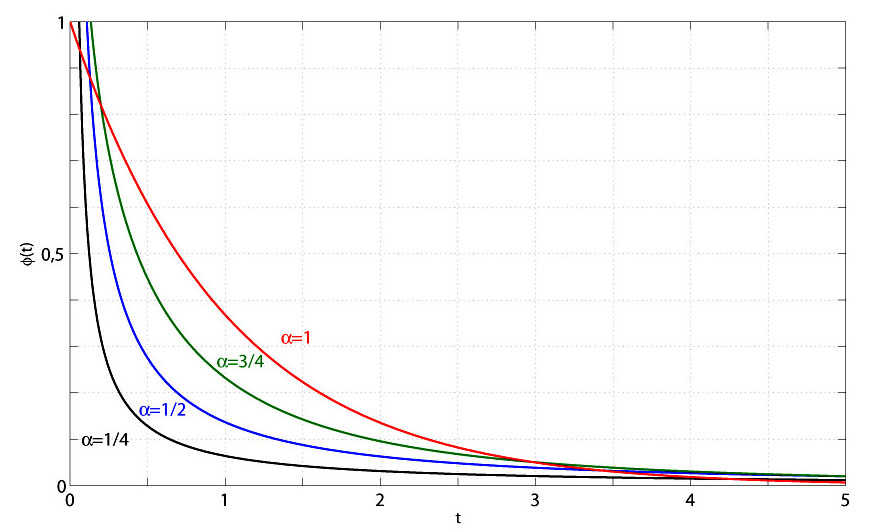}
\end{center}
%\vskip -0.5truecm
 %\caption
 \caption{Plots  of  %% Mittag-Leffler type 
  $\phi _\alpha (t)$  with $\alpha = 1/4, 1/2, 3/4, 1$ versus $t$;
 for $0\le t \le 5$.}
 % \label{fig:1.2}
 \label{fig5}
 \end{figure}

%\vsp
As a matter of fact
a straightforward generalization of the classical  Mittag-Leffler function
is obtained  by replacing the additive constant $1$
in the argument of the Gamma function in (\ref{2.1}) by an
arbitrary complex parameter $\beta  \,. $
It was  formerly considered in 1905 by Reference \cite{Wiman 05a}
and soon later by  Mittag-leffler  himself, almost incidentally  in one of his notes. 
 Later,  in the 1950's,  such generalization  was investigated by Humbert  and Agarwal,
with respect to the  Laplace transformation, see  
References \cite{Humbert 53}, \cite{Agarwal 53}, \cite{Humbert-Agarwal 53}.
Usually,  when   dealing with Laplace transform pairs,
the parameter $\beta  $ is required to be real and positive  like  $\alpha$.
\vsp
For this function we agree to use the  notation
\bee %% $$
E_{\alpha, \beta  } (z ) := \sum_{n=0}^\infty
  \frac{z^n}{ \Gamma(\alpha n+ \beta   )}\,,
\quad \Re \alpha >0\,,\;\, \beta  \in \CC\,,
   \quad z \in\CC
  \,. 
 \label{E.22} %%  \eqno(E.22) $$
 \ee
Of course $E_{\alpha,1}(z) \equiv E_\alpha(z)$.
The series is still convergent for all the complex plane $\CC$
so the function (\ref{E.22}) is still entire for $\Re (\alpha) >0 $ for any
$\beta \in \CC$ with order $1/{\Re \alpha}$  so   the additional parameter %$\beta$ does not 
play any role on this respect.
However the Laplace transform pairs concerning  the Mittag-Leffler function (\ref{E.22})
and its derivative are known to be 
with $\alpha, \beta >0$ and  $\Re (s) > |\lambda|^{1\alpha}$, 
see, for example, \cite{Podlubny BOOK99}, \cite{Mainardi BOOK10}, 
\cite{GKMR BOOK20},
\bee  %% $$
%%e_{\alpha,\beta }  (t;\lambda):=
   t^{\beta   -1}\,  E_{\alpha,\beta  } \left(-\lambda \, t^\alpha\right)
   \, \div\, 
 \frac{s^{\alpha -\beta  }}{ s^\alpha +\lambda}
 = \frac{s^{-\beta}} { 1 + \lambda s^{-\alpha}}
           \,. 
\label{E.53}   %%           \eqno(E.53)$$
 \ee
 and 
\bee %%  $$ 
t^{\alpha k + \beta -1}\, E^{(k)}_{\alpha, \beta}( \lambda t^\alpha)
\, \div \, 
\frac{ k!\, s^{\alpha-\beta}}{(s^\alpha - \lambda )^{k+1}}\,,
\q k = 0,1, 2, \dots  \,. 
\label{E.70} %%  \eqno(E.70)$$
\ee
We also note the following relation concerning the first derivative of the classical 
Mittag-Leffler function with the two-parameter Mittag-Leffler function 
usually overlooked by several authors but of easy proof:
\bee %%$$ 
 \phi _\alpha (t) := {\ds t^{-(1-\alpha  )}\, E_{\alpha  ,\alpha } \left(-  t^{\alpha  }\right)}
=  - {\ds  \frac{d}{dt} E_\alpha   \left (- t^{\alpha}\right)}\,,
\quad t\ge 0,
\quad  0<\alpha  <1\,.
\label{1.45}  %%\eqno(1.45)$$
 \ee
 We report the plot of the function $\phi_\alpha(t)$  herewith in Figure~\ref{fig5}. 
 
% \vsp
We note that Mittag-Leffler functions with more than two parameters  were
 also dealt by several authors as pointed out in \cite{GKMR BOOK20}.
  In particular, for the 3-parameter Mittag-Leffler function (known as Prabhakar 
function) and related operators we refer the reader to the 
recent survey by Giusti et al.~\cite{Giusti FCAA20}
and references therein. 
Kiryakova  has dealt 
in a number of papers the multi-index Mittag-Leffler functions, see for example  
\cite{Kiryakova CAMWA10}.    

\section{The fractional Poisson process and the Mittag-Leffler function}

Hereafter we describe how the Mittag-Leffler function enters in the so-called fractional Poisson process.
We are following the original approach by Mainardi et al. in 
\cite{Mainardi-Gorenflo-Scalas JV04} where the fractional Poisson process 
is referred to as the renewal process of the Mittag-Leffler type.
However, an independent approach to the fractional Poisson process was given
for example,  by 
%Uchaikin \cite{Ucaikin FPP} and 
 Laskin in ~\cite{Laskin FPP03}.

\subsection{Essentials of Renewal Theory}

The concept of  {\it renewal process} has been developed as a
stochastic model for describing the class of counting processes
for which the times between successive  events are
independent  identically distributed  ($iid$)
non-negative random variables, obeying a given probability law.
These times are  referred to as waiting times
or inter-arrival times.
%% they  will be here denoted by $T_1, T_2, \dots$
%% These processes  generalize the classical Poisson process, that
%% is known to be   a counting process
%% where the waiting times are exponentially distributed.
%% and when it fails we immediate replace it with a new one.
For more details see, for example, the classical
treatises by
%Khintchine \cite{Khintchine QUEUEING60},
Cox \cite{Cox RENEWAL67},
%Gnedenko \& Kovalenko \cite{Gnedenko-Kovalenko QUEUEING68},
Feller \cite{Feller BOOK71}.
% and the  recent book by Ross \cite{Ross PROBMOD97}.
%% More recent books that treat in great detail such processes include
%% \cite{Ross STOCHPROC96,Ross PROBMOD97} and
%% \cite{BeicheltFatti STOCHPROC02}.
%%%%%%%%%%%%%\vfill\eject%%%%%%%%%%%%%%%

For a renewal process having waiting times $T_1,T_2, \dots$, let
\bee %% $$
t_0 = 0\,, \q   t_k= \sum_{j=1}^k T_j\,, \q k \ge 1\,. 
\label{1.1V} %% \eqno(1.1) $$
\ee
That is $t_1 =T_1$ is the time of the first renewal,
$t_2 = T_1 +T_2$ is the time of the second renewal and so on.
In general $t_k$ denotes the $k$th renewal.

The process is specified if we know the probability law for the waiting
times. In this respect we introduce the
{\it probability density function} ($pdf$)
$\phi(t)$
 and  the   (cumulative) distribution function $\Phi(t)$
so defined:
\bee %%$$
 \phi (t) := \frac{d}{dt} \Phi(t) \,,\q
   \Phi(t) := P \lt( T \le t\rt) = \int_0^t \phi (t')\, dt'\,.
\label{1.2V} %%\eqno(1.2)$$
\ee
 When the non-negative random variable  represents
 the lifetime of technical systems, it is common to refer
to $\Phi(t)$ as to the {\it failure probability}
and to
\bee %% $$ 
\Psi(t) := P \lt(T > t\rt) = \int_t^\infty \phi (t')\, dt' = 1-\Phi(t)\,,
\label{1.3V}  %% \eqno(1.3)$$
\ee
as to the {\it survival probability}, because $\Phi(t)$ and $\Psi(t)$ are
the respective probabilities that the system does or does not fail
in $(0, T]$. %%  see \eg \cite{BeicheltFatti STOCHPROC02}.
%%%%
A relevant quantity is the {\it counting function} $N(t)$
defined as
\bee %% $$ 
N(t):= \hbox{max} \lt\{k | t_k \le t, \;k = 0, 1, 2, \dots\rt \}\,,
\label{1.4V}  %%\eqno(1.4)$$
\ee
that represents the effective number of events before or at instant $t$.
In particular we have $\Psi(t) = P \lt(N(t) =0\rt)\,.$
%%%%%%%%%
Continuing in the general theory we set
$F_1(t) = \Phi(t)$, $f_1(t) = \phi (t)$, and in general
\bee %%$$ 
F_k(t) :=  P\lt(t_k = T_1+ \dots +T_k \le t \rt)\,,
 \; f_k(t) = \frac{d}{dt} F_k(t)\,, \; k\ge 1\,, 
 \label{1.5V} %%\eqno(1.5)$$
 \ee
thus  $F_k(t)$ represents the probability that the sum of the first $k$
waiting times is less or equal  $t$
and $f_k(t)$ its density. Then, for any fixed $k\ge 1$  the normalization
condition  for $F_k(t)$ is fulfilled
%% must hold  as $t \to \infty$, namely
because
\bee %% $$
         \lim_{t \to \infty} F_k(t) =
  P\lt(t_k = T_1+ \dots +T_k < \infty \rt)= 1 \,.
\label{1.6V} %%                 \eqno(1.6)$$
\ee
In fact,    the sum of $k$ random variables each of which is
    finite with probability 1 is finite with probability 1 itself.
%%%%%
By setting for consistency $F_0(t) \equiv 1$ and $f_0(t) = \delta(t)$, where
for the Dirac delta  generalized function
in $\RR^+$ we assume the {\it formal representation}
$$\delta(t) := \frac{t^{-1}}{\Gamma(0)}\,, \q t\ge 0\,,$$
we also note that for $k \ge 0$ we have
\bee %% $$ 
P\lt(N(t) =k\rt) := P\lt(t_k \le t \,,\,t_{k+1} > t\rt)
=  \int_0^t f_k(t')\, \Psi (t-t')\, dt'\,.
%%  = f_k(t) \, * \, \Psi(t) \,.
\label{1.7V} %%\eqno(1.7)$$
\ee
%% where $\Psi(t)$ is the survival probability defined in (1.3).
%%%%%%%%%%%%%%%%%\vfill\eject%%%%%%%%%%%

We now find it convenient to introduce the simplified $\, *\,$ notation
for the Laplace convolution between two causal well-behaved
(generalized) functions $f(t)$ and $g(t)$
$$   \int_0^t   f(t')\, g(t-t')\, dt' = \lt(f \,*\, g\rt) (t) =
     \lt(g \,* \, f\rt) (t)  = \int_0^t   f(t-t')\, g(t')\, dt'\,. $$
Being  $f_k(t)$ the $pdf$ of the sum of the  $k$
$iid$ random variables
$T_1,  \dots, T_k$
with $pdf$ $\phi (t)\,, $ we easily recognize that
$f_k(t)$ turns out to be the $k$-fold convolution of $\phi(t)$
with itself,
\bee %% $$ 
f_k(t) =  \lt(\phi^{*k}\rt) (t)\,, 
\label{1.8V} %%\eqno(1.8)$$
\ee
so Equation (\ref{1.7V})  simply reads:
\bee %% $$   
 P\lt(N(t) =k\rt) = \lt(\phi^{*k} \,*\, \Psi\rt)(t)\,. 
 \label{1.9V} %%\eqno(1.9)$$
 \ee
%%%%%%%%%%%%%
Because of the presence of Laplace convolutions
a renewal process is suited for the Laplace transform method.
Throughout this paper we will denote
by   $\widetilde f(s)$
the Laplace transform
of a sufficiently well-behaved (generalized) function
$f(t)$  according to
$$   {\L} \lt\{ f(t);s\rt\}=   \widetilde f(s)
 = \int_0^{+\infty} \e^{\ds \, -st}\, f(t)\, dt\,,
\q s > s_0\,,
$$
and for $\delta (t)$
consistently we will have
$  \widetilde \delta (s) \equiv 1\,. $
{\it} Note that for our purposes we agree to take $s$ real.
%% in the Laplace transform.
%% TO BE CONTINUED
We recognize  that
(\ref{1.9V}) reads in the Laplace domain
\bee %% $$ 
   \L \{P\lt( N(t)=k\rt); s\}
= \lt[{\widetilde \phi (s)} \rt]^k \,   \widetilde \Psi (s)
             \,,
\label{1.10V} %%       \eqno(1.10)$$
\ee
where, using (\ref{1.3V}),
\bee %%$$ 
\widetilde \Psi (s)  =
\frac{ 1- \widetilde \phi (s)} {s}\,.
\label{1.11V} %%\eqno(1.11)$$
\ee
%%%%%%%%%\vfill\eject%%%%%%%%%%

%%%%%%%

\subsection{The Classical Poisson Process as a Renewal Process}

The most celebrated   %%% simple
renewal process is the Poisson process
characterized by a waiting time $pdf$ of exponential type,
\bee %% $$ 
 \phi (t) = \lambda \, \e^{-\lambda t}\,,\q \lambda >0\,, \q t\ge 0\,.
\label{2.1V}  %%\eqno(2.1)$$
\ee
The process has {\it no memory}.
Its  moments turn out to be
 \bee %% $$ 
 \langle T\rangle = \frac{1}{\lambda }\,, \q
   \langle T^2 \rangle = \frac{1}{\lambda^2 }\,,
\q \dots \,,  \q \langle T^n \rangle = \frac{1}{\lambda^n }\,, \q \dots \,,
\label{2.2V}  %%\eqno(2.2)$$
\ee
and the {\it survival probability} is
\bee %%   $$
 \Psi(t) := P\lt(T>t\rt) = \e^{-\lambda t}\,, \q t \ge 0 \,.
\label{2.3V}  %% \eqno(2.3)$$
\ee
%% The exponential distribution is characteristic of a process without memory.
%% In fact
%% $$ \qq P \l( T > s+t | T>s \r) =
%% \frac{P \l( T > s+t \,,\, T>s \r)} {P \l( T>s \r)}
%% = \frac{P \l( T > s+t\r) } {P \l( T>s \r)} \qq \eqno(2.4)$$
%% $$ = \frac{\e^{-\lambda (s+t)}}{ \e^{-\lambda s}} =
%%    \e^{-\lambda t} = P\l(T >t\r)\,. $$
We know that the probability  that $k$ events occur in the
interval of length $t $ is
\bee %% $$ 
P\lt( N(t) = k\rt) =   \frac{(\lambda t)^k}{k!} \, \e^{-\lambda t}
\,, \q t \ge 0\,, \q k = 0,1, 2, \dots\;.
\label{2.4V} %% \eqno(2.4)$$
\ee
%%% Furthermore the Poisson process turns out to be a renewal process
%%% with linear renewal function, namely
%% $$ m(t) = \lambda \, t\,,  \q t \ge 0
%% \,.$$
The probability distribution related to the sum of $k$ $iid$
exponential random variables  is known to be
the so-called {\it Erlang distribution}  (of order $k$).
The corresponding density (the {\it Erlang} $pdf$) is thus
\bee %% $$ 
f_k(t) = \lambda \,  \frac{(\lambda t)^{k-1}}{(k-1)!}
    \, \e^{-\lambda t}\,,\q t \ge 0 \,,
\q k =1,2, \dots \,, 
\label{2.5V}  %%\eqno(2.5)$$
\ee
so that the  Erlang distribution function of order $k$ turns out
to be
\bee %% $$ 
F_k(t) = \int _0^t f_k(t')\, dt' =
     1 -  \sum_{n=0}^{k-1}
      \frac{(\lambda t)^n}{n!} \, \e^{-\lambda t}  =
 \sum_{n=k}^{\infty}
      \frac{(\lambda t)^n}{n!} \, \e^{-\lambda t}\,,\q t \ge 0
\,. 
\label{2.6V} %%\eqno(2.6)$$
\ee
In the limiting case $k=0$ we recover
%%% (for consistency)
$f_0(t) = \delta(t),\; F_0(t) \equiv 1,\; t\ge 0$.
%% We note hat $F_k(t)$ represents
%%%% the probability that $k$ or more events occur in $(0,t]$, namely
%% the probability that the sum of the first $k$
%% waiting times are less or equal  $t$.

The results  (\ref{2.4V})--({\ref{2.6V}) can easily obtained
by using the technique of the Laplace transform sketched in
the previous section
noting that for the Poisson process we have:
\bee %% $$ 
\widetilde\phi(s) =
\frac{\lambda } {\lambda +s}\,,
 \q  \widetilde\Psi(s) =  \frac  {1 } {\lambda +s}\,,  
 \label{2.7V} %% \eqno(2.7)$$
 \ee
and for the Erlang distribution:
\bee %% $$
  \widetilde f_k(s) =  [\widetilde\phi(s)]^k =
    \frac  {\lambda^k } {(\lambda +s)^k}\,,
\q    \widetilde F_k(s) =  \frac{[\widetilde\phi(s)]^k}{s} =
    \frac  {\lambda^k } {s (\lambda +s)^k}\,.
\label{2.8V}  %%\eqno(2.8)$$
\ee

We also recall that the survival probability
for the Poisson renewal process
obeys the ordinary differential equation  (of relaxation type)
\bee %% $$
     \frac{d}{dt} \Psi(t) = -\lambda \Psi(t)\,, \q t \ge 0\,; \q
\Psi(0^+) =1\,. 
\label{2.9V}  %%\eqno(2.9)$$
\ee

 \subsection{The Renewal  Process of Mittag-Leffler Type}

A ``fractional'' generalization of the   Poisson renewal process
is simply obtained by generalizing the differential Equation  (\ref{2.9V})
replacing there the first derivative with the integro-differential operator
$\, _*D_t^\beta$ that  is interpreted  as
the fractional derivative
of order $\beta $ in  Caputo's sense, see Section 2.
We write, taking for simplicity $\lambda =1$,
\bee %%$$ 
      \, _*D_t^\beta \, \Psi(t) =
- \Psi(t)\,, \q t >0\,,\q 0<\beta \le 1\,; \q \Psi(0^+) =1\,. 
\label{3.1V}  %%\eqno(3.1)$$
\ee
We also allow
the limiting case $\beta =1$ where all the results of the previous
section (with $\lambda =1$) are expected to be recovered.

%\vsp
For our purpose we need to recall the Mittag-Leffler function
as the natural ``fractional'' generalization
of the exponential function, that characterizes the Poisson process.
We again recall that the Mittag-Leffler function of parameter  $\beta\,$
is defined in the complex plane by the power series
\bee %% $$
 E_\beta (z) :=
    \sum_{n=0}^{\infty}\,
   {z^{n}\over\Gamma(\beta\,n+1)}\,, \q \beta >0\,, \q z \in \CC\,,
\label{3.2V} %%  \eqno  (3.2)$$
\ee
as stated in Section 2 where the parameter was denoted by $\alpha$.
%%%%

The solution  of Equation (\ref{3.1V}) is known to be, see 
Section 3 
%\eg
%\cite{CaputoMainardi RNC71,GorMai CISM97,Mainardi CHAOS96},
\bee %% $$
\Psi(t) =  E_\beta (-t^\beta)\,, \q t \ge 0\,, \q 0<\beta \le 1\,,
\label{3.3V} %%\eqno (3.3)$$ 
\ee 
so
\bee %% $$
 \phi(t) :=    -   \frac{d} {dt} \Psi(t) =
            -   \frac{d}{dt}  E_\beta (-t^\beta) \,, \q t \ge 0
 \,, \q 0<\beta \le 1\,.
 \label{3.4V} %% \eqno (3.4)$$
 \ee
Then, the corresponding Laplace transforms read
\bee %% $$ 
\widetilde \Psi(s) =
 {s^{\beta-1} \over 1+ s^\beta}\,,   \q
 \widetilde \phi(s)= {1 \over 1+  s^{\beta}}\,,\q
    0<\beta \le 1\,.
   \label{3.5V} %% \eqno(3.5) $$
   \ee
%%%%%%%%%
Hereafter, we find it convenient to summarize
the  most relevant features  of the functions $\Psi(t)$ and $\phi(t)$
when  $0< \beta <1\,.$
We begin  to quote their series expansions convergent in all of $\RR$
suitable for   $t \to 0^+ $
and  their asymptotic representations  for  $t\to \infty $,
\bee %% $$ 
\Psi(t)
   = {\ds \sum_{n=0}^{\infty}}\,
  (-1)^n {\ds {t^{\beta n}\over\Gamma(\beta\,n+1)}}
 \,\sim \,  {\ds {\sin \,(\beta \pi)\over \pi}}
  \,{\ds  {\Gamma(\beta)\over t^\beta}}\,,
 %%%% & $\, t\to \infty \,,$} \q 0<\beta <1\,,
\label{3.6V} %%     \eqno(3.6) $$
\ee
and
\bee %% $$ 
\phi(t)
= {\ds {1\over t^{1-\beta}}}\, {\ds \sum_{n=0}^{\infty}}\,
  (-1)^n {\ds {t^{\beta n}\over\Gamma(\beta\,n+\beta )}}
 \, \sim \,  {\ds {\sin \,(\beta \pi)\over \pi}}
  \,{\ds  {\Gamma(\beta+1)\over t^{\beta+1}}}\,.
%%%  & $\, t\to \infty \,.$} \q 0<\beta <1\,,
\label{3.7V} %%     \eqno(3.7) $$
\ee
%%%%%
In contrast %%%% Contrarily %% At variance 
to the Poissonian case  $\beta=1$,
in the case  $0<\beta <1$ for large $t$
the functions $\Psi(t)$ and $\phi(t)$
no longer decay   exponentially  but algebraically.
%% %%%({\it power-law decay}).
As a consequence of the power-law asymptotics
%% for the waiting time $pdf$
the process turns be no longer Markovian but of long-memory type.
However, we recognize that for $0<\beta <1$ both  functions
 $\Psi(t)$, $\phi(t)$
%% even if they lost their exponential decay by exhibiting power-law tails
%% for large times,
keep   the ``completely monotonic'' character of the Poissonian case.
%% of the exponential function of $\exp (-t)\,,$
as can be simply derived from Section 2.
We  recall that {\it complete monotonicity}  of  our   functions
 $\Psi(t)$ and $\phi(t)$  means
 \bee %% $$ 
 (-1)^n {d^n\over dt^n}\, \Psi  (t) \ge 0\,,  \q
   (-1)^n {d^n\over dt^n}\, \phi  (t) \ge 0\,,
\q n=0,1,2,\dots   \,, \q t \ge 0\,,     
\label{3.8V}  %%\eqno(3.8)$$
\ee
or equivalently, their representability as real Laplace transforms
of non-negative   generalized
functions (or measures).
%%  see \eg \cite{GorMai  CISM97}. %%,Mainardi BONN00}.
%%%%%%%%%%

For the generalizations
of Equations (\ref{2.4V})--(\ref{2.6V}),  characteristic
of the Poisson and Erlang distributions respectively,
we must point out the   Laplace transform pair
\bee 
 t^ {\beta \,k}\, E_\beta ^{(k)}
  (-t^\beta )  \, \div\, 
        \frac{ k!\, s^{\beta -1}}{(1+s^\beta )^{k+1}}
\,, \q \beta >0 \,, \q k = 0,1, 2, \dots \,,
\label{3.9V} %% \eqno(3.9)$$
\ee
with $ {\ds E_\beta ^{(k)}(z) := \frac{d^k}{dz^k}  E_\beta(z)}\,, $
that can be deduced  from the book by Podlubny,
see Equation (1.80) in  Reference \cite{Podlubny BOOK99}.
%%%%%%%\vfill\eject%%%%%%%%%%%%%
Then, by using the Laplace transform pairs  (\ref{3.5}) and
Equations (\ref{3.3V}), (\ref{3.4V}), (\ref{3.9V})
in Equations  (\ref{1.8V}) and (\ref{1.9V}),
 we have  the {\it generalized Poisson distribution},
\bee %% $$ 
P\lt( N(t) = k\rt) =   \frac{ t^{ k\, \beta}}{k!} \,
  E_\beta^{(k)} (- t^\beta)
\,, \q k = 0, 1, 2, \dots 
\label{3.10V} %%\eqno(3.10)$$
\ee
and the {\it generalized Erlang} $pdf$'s   (of order $k \ge 1$),
\bee %% $$
 f_k(t) = \beta \,  \frac{ t^{k\beta-1}}{(k-1)!}
    \, E_\beta^{(k)} (- t^\beta)   %% \,, \q k =1,2, \dots
\,. 
\label{3.11V} %%\eqno(3.11)$$
\ee
The {\it generalized  Erlang distribution functions} turn out
to be
\bee %% $$ 
F_k(t) = \int _0^t f_k(t')\, dt' =
     1 -  \sum_{n=0}^{k-1}
      \frac{ t^{n \beta}}{n!} \, E_\beta^{(n)} (- t^\beta)  =
 \sum_{n=k}^{\infty}
      \frac{t^{n\beta}}{n!} \, E_\beta^{(n)} (- t^\beta)
\,. 
\label{3.12V}%% \eqno(3.12)$$
\ee

%% It may be instructive to consider the special case $\beta =1/2$
%% for which it is known that
%% $$ E_{1/2} (-\sqrt{t}) =
%%    \e^{\ds \, t}\, \hbox{erfc} (\sqrt{t}) =
%% \e^{\, \ds t}\, {2\over \sqrt{\pi}}\,
%%  \int_{\sqrt {t}}^\infty \e^{\, \ds -u^2}\,du  \,,\q t\ge 0\,,
%% \eqno(3.13)$$
%% where $ \, \hbox{erfc}\,$ denotes the {\it complementary error} function.
%% In this case we can take profit of the recurrence relations for
%% repeated integrals of the error functions,
%% see \eg \cite{AS 65},\S 7.2, pp 299-300,
%% to compute the derivatives of the Mittag-Leffler.
%% For this purpose we recall:
%% $$ \frac{d^n}{dz^n} \l( \e^{z^2}\, \erfc (z) \r) =
%%   (-1)^n \, 2^n\, n!\, \e^{z^2}\, I^n \, \erfc (z)\,,\eqno(3.14)$$
%% where
%% $$I^n \, \erfc (z) =\int_z^\infty I^{n-1} \,\erfc (\zeta)\, d\zeta\,,
%% \q n= 0, 1, 2, \dots \eqno(3.15)$$
%% with
%% $$ I^{-1} \, \erfc (z) = \frac{2}{\sqrt {\pi}} \, \e^{-z^2}\,, \q
%%   I^{0} \, \erfc (z) = \erfc (z)\,. \eqno(3.16)$$
%%%%%%%%%%%%%%%%\vskip 0.5 truecm %%%%%%%%%%%

\section{The Gnedenko-Kovalenko  theory of thinning
and the Mittag-Leffler function}
%%% SECTION 2 in BH Bad-Honnef chapter

The {\it thinning} theory for a renewal process has been considered
in detail by  Gnedenko and Kovalenko \cite{Gnedenko-Kovalenko 68}
in the first edition of their book on Queue theory of 1968.
%\cite{Gnedenko-Kovalenko 68}.
However, the connection with the Laplace transform of the  Mittag-Leffler 
function outlined at the end of this section in Equations (\ref{2.9BH}),
(\ref{2.10BH}), 
see also \cite{Gorenflo-Mainardi Bad-Honnef06}
and  \cite{Gorenflo-Mainardi TIMS},
is surprisingly not present in the second edition of 
the book by Gnedenko \& Kovalenko in 1989.

%\vsp
We must note that other authors, like Sz\'antai
\cite{Szantai 71a,Szantai 71b} speak of {\it rarefaction}
in place of thinning.

%\vsp
Let us sketch here the essentials of this theory:
in the interest of transparency and easy
readability we avoid the possible decoration
of the relevant power law by multiplying it with a 
{\it slowly varying function}.
%% {\bf Remark}
%%As usual  we call a (measurable) positive function $a(y)$ 
%%{\it slowly varying at zero} if
%% $a(cy)/a(y) \to 1$ with $y \to 0^+$ for every $c>0$,
%% {\it slowly varying at infinity}
%% if $a(cy)/a(y) \to 1$ with $y \to +\infty$ for every $c>0$.
%%A standard example of a slowly varying function at zero and at infinity is 
%% $|\log y|^\gamma$, with $\gamma \in \RR$.

%\vsp
Denoting by $t_n$,\, $n=1,2,3, \dots$
the time instants of events of a renewal process, assuming
$0=t_0<t_1<t_2<t_3<\dots $,
with $i.i.d.$  waiting times $T_1 = t_1\,,\,T_k = t_{k}-t_{k-1}$ for $k\ge 2$,
%%  given as $iid$ random variables
(generically denoted by T),
{\it thinning} (or {\it rarefaction})
means that for each positive
index  $k$  a decision is made:
the event happening in the instant $t_k$ is deleted with probability $p$
or it is maintained with probability $q=1-p$,  $0<q<1$.
This procedure produces a {\it thinned} or {\it rarefied} renewal process
with fewer events (very few events
if $q$ is near zero, the case of particular interest)
in a moderate span of time.

%\vsp
To compensate for this loss  we change the unit of time
so that we still have not very few  but still a moderate number
of events in a moderate span of time.
Such change of the unit of time is
equivalent to  rescaling the waiting time,
multiplying it with a positive factor $\tau $ so that we have
waiting times $\tau T_1,\tau T_2, \tau T_3, \dots$, and
instants $\tau t_1,\tau t_2, \tau t_3,\dots$, in the
rescaled  process.
Our intention is, vaguely speaking, to dispose
on $\tau $ in relation to the rarefaction parameter $q$
in such a way that for $q$ near zero in some sense
the ``average'' number of events per unit of time
remains unchanged. In an asymptotic sense
we will make these considerations precise.

%\vsp
Denoting by $F(t) = P(T\le t)$ the probability distribution function
of the (original) waiting time $T$,   by $f(t)$ its density
($f(t)$ is a generalized function generating a probability measure)
so that
$F(t) = \int_0^t f(t') \, dt'$,  and analogously by
$F_k(t)$  and $f_k$(t) the distribution %% probability distribution function
and density, %%% the probability density function,
respectively,  of the sum of $k$  waiting times, we have recursively
\bee %%$$ f_1(t) = f(t) \,,\q
      f_k(t) = \int_0^t  f_{k-1}(t-t') \, dF(t') \,,\;
\hbox{for} \; k \ge 2\,. 
\label{2.1BH} %% \eqno(2.1)$$
\ee
Observing that after a maintained event the next one of the
original process is kept with probability $q$ but dropped in favour
of the second-next with probability $p\,q$
and, generally, $n-1$ events are dropped  in favour of the
$n$-th-next with probability $p^{n-1}\,q$,
we get for   the waiting time density of the thinned process
the formula
\bee %% $$
 g_q(t) = \sum_{n=1}^\infty q\, p^{n-1}\, f_n(t)\,.
\label{2.2BH} %% \eqno(2.2)$$
\ee
%%%%%%%%
With the modified waiting time $\tau \,T$ we have
$$ P(\tau T\le t) = P(T\le t/\tau ) = F(t/\tau )\,, $$
hence the density $f(t/\tau )/\tau $, and analogously for the
density of the sum of $n$ waiting times
$f_n(t/\tau )/\tau $.
The density of the waiting time of the rescaled (and thinned) process
now turns out as
\bee %%$$ 
g_{q,\tau} (t) = \sum_{n=1}^\infty q\, p^{n-1}\, f_n(t/\tau)/\tau \,.
\label{2.3BH} %%\eqno(2.3)$$
\ee
\vsp
In the Laplace domain we have
 $\widetilde f_n(s) = \left(\widetilde f(s)\right)^n\,,$
hence (using $p =1-q$)
\bee %% $$ 
\widetilde g_q (s)=
\sum_{n=1}^\infty q\, p^{n-1}\,\left(\widetilde f(s)\right)^n
 = \frac{ q\, \widetilde f(s)}{ 1 - (1-q)\, \widetilde f(s)}
 \,,
 \label{2.4BH} %%\eqno(2.4)$$
 \ee
from which by Laplace inversion we can, in principle, construct
the waiting time density of the thinned process.
By  rescaling we get
\bee %% $$
\widetilde g_{q,\tau}(s)=
\sum_{n=1}^\infty q\, p^{n-1}\,\left(\widetilde f(\tau s)\right)^n
 = \frac{ q\, \widetilde f(\tau s)}{ 1 - (1-q)\, \widetilde f(\tau s)}
 \,.
\label{2.5BH} %% \eqno(2.5)$$
\ee
\vs
\noindent
% \subsection*{Infinite thinning under proper rescaling}
%\vskip0.25truecm \noindent
     Being interested in stronger and  stronger
thinning ({\it infinite thinning})
let us now  consider 
 a scale of processes with  the parameters $\tau  $ (of {\it rescaling})
and $q$ (of {\it thinning}), with $q$ tending
to zero  {\it under a scaling relation $q = q(\tau) $
  yet to be specified}.

%\vsp
%% To get handsome results we assume that the original waiting time
%% either has a finite expectation or is distributed according to a power law.
We have essentially two cases for the waiting time distribution:
its expectation value is finite or infinite.
In the first case we put
\bee %% $$ 
\lambda  = \int_0^\infty t'\,  f(t')\, dt' < \infty \,.
\label{2.6aBH} %%\eqno (2.6a)$$
\ee
In the second case we assume a queue of power law type
(dispensing with a possible decoration by 
a function slowly varying  at infinity)
\bee %% $$
  \Psi(t) :=  \int_t^\infty f(t')\, dt'
\sim \frac{c}{\beta } t^{-\beta }\,, \; t\to \infty \quad
\hbox{if}\quad 0<\beta <1\,.
\label{2.6bBH} %%{\eqno(2.6b)$$
\ee
Then, by the Karamata theory (see References \cite{Feller BOOK71,Widder BOOK46})
the above  conditions mean in the Laplace domain
\bee %% $$
\widetilde f(s) =1-\lambda  \, s^\beta  + o\left( s^\beta \right)\,,
 \q \hbox{for} \q s \to 0^+\,,
 \label{2.7BH} %%\eqno(2.7)$$
 \ee
with  a positive coefficient $\lambda $ and $0<\beta \le 1$.
The case $\beta =1$  obviously corresponds to the situation with finite
first moment (2.6a), whereas the case $0<\beta <1$ is related
to a power law   queue with
$c= \lambda \,\Gamma(\beta +1)\,\sin(\beta \pi)/\pi\,.$

%\vsp
Now, passing  to the limit of $q \to 0$ of infinite thinning under the scaling relation
\bee %% $$ 
 q = \lambda \, \tau ^\beta \,, \q 0<\beta \le 1\,, 
 \label{2.8BH} %% \eqno(2.8)$$
 \ee
between the positive parameters $q$ and $\tau $,
the Laplace transform of the rescaled density
$\widetilde {g_{q,\tau }}(s)$ in (\ref{2.5BH})
 of the thinned process tends for fixed $s$ to
\bee %% $$  
\widetilde g(s) = \frac{1}{1+s^\beta}\,, 
\label{2.9BH} %% \eqno(2.9)$$
\ee
which corresponds to the Mittag-Leffler density
\bee %%$$ 
g(t) = - \frac{d}{dt} E_\beta (-t^\beta )
= \phi^{ML}(t)
\,. 
\label{2.10BH} %%\eqno(2.10)$$
\ee
Let us remark that Gnedenko and Kovalenko obtained
(\ref{2.9BH})  as the Laplace transform of the limiting density
but did not identify it as the Laplace transform of a
Mittag-Leffler type function.
Observe that in the special case $\lambda < \infty$ 
we have $\beta=1$,
hence as the limiting process the Poisson process, as formerly shown
in 1956 by R{\'e}nyi  \cite{Renyi 56}.

%%%%%%%%%%%\vfill\eject

%%%%%%%%%

\section{The Continuous Time Random Walk (CTRW) and the Mittag-Leffler function}  %%% SECTION 3

The name    {\it continuous time random walk} (CTRW)  became
popular in physics after Montroll and Weiss  (just to cite  the pioneers)
published a celebrated series
of papers on random walks for modelling
diffusion processes on lattices, see, for example,
Reference \cite{Montroll-Weiss 65}, and
the book by Weiss \cite{Weiss BOOK94} with  references therein.
CTRWs are rather good and general phenomenological models for diffusion,
including processes of anomalous transport,
that can be understood  in the framework of
the classical renewal theory.
%% as stated \eg in  the booklet by Cox  \cite{Cox_RENEWAL67}.
In fact a CTRW can be considered
as a compound renewal process (a simple renewal process with reward) or
 a  random walk {\it subordinated}
to a simple  renewal process.
Hereafter we will mainly follow the approach by Gorenflo \& Mainardi, 
see, for example,  Reference \cite{Gorenflo-Mainardi KLM12}.
%%\vsp
%% Basic notions of the CTRW theory,  that hereafter we briefly recall
%% for the readers' convenience, are the master  equation (in integral form)
% for the sojourn probability density,
%% its Fourier-Laplace representation (known as the Montroll-Weiss formula)
%% and its series representation.  %% (known as the Cox formula)

%\vsp
A spatially one-dimensional CTRW  %% continuous time random walk (CTRW)
is generated by a sequence
of  independent identically  distributed ($iid$)
 positive  random  waiting times $T_1, T_2, T_3, \dots ,$
each having the same probability density function
 $\phi(t)\,,$   $\, t>0\,, $ and
a sequence of $iid$ random jumps $X_1, X_2, X_3, \dots, $
in $\RR\,,$ each having the same probability density
$w(x)\,,$ $\, x\in \RR\,.$

%\vsp
Let us remark that, for ease of language, we use
the word density also for generalized functions
in the sense of Gel'fand \& Shilov \cite{GelfandShilov 64},
that can be interpreted as probability measures.
Usually the {\it probability density functions} are abbreviated
by  $pdf$.
We recall that $\phi (t) \ge 0$ with $\int_0^\infty \phi (t)\, dt =1$
and  $w(x)  \ge 0$ with $\int_{-\infty}^{+\infty} w (x)\, dx =1$.
%%%%%%%%%%%% \newpage

%\vsp
%%%%%%%%%%%%
Setting
$t_0=0\,,$ $\, t_n = T_1+T_2 + \dots T_n$ for $n \in \NN\,,$
%%% $0< t_1<t_2 < \dots\,,$
the wandering particle %%% (or the wanderers?)
%%  starts at point $x=0$ in instant $t=0$ and
makes a jump
of length $X_n$ in instant $t_n$,
so that its position is $x_0=0$ for
$0\le t <T_1= t_1\,,$       and
$x_n =    X_1 + X_2 + \dots X_n\,,$
for $ t_n \le t < t_{n+1}\,. $
%% An often (also by us) required assumption  is that
We require the distribution of
the waiting times
and  that of the jumps to be independent
of each other.
So, we have a compound renewal process (a renewal process with
       reward), compare  Reference \cite{Cox RENEWAL67}.

%\vsp
%%%%%%%%%%%%%%%%%%%
By natural probabilistic arguments we arrive at the
{\it integral  equation} for the probability density   $p(x,t)$
(a density with respect to the variable $x$)
of the particle being in point $x$ at instant $t\,, $
%see  \eg
%\cite{Gorenflo-%Mainardi_INDIA03,Gorenflo_KONSTANZ01,Mainardi_PhysicaA00,%
% Scalas_PhysA05,SGM_00,Scalas_PRE04},
\bee %%$$
   p(x,t) =  \delta (x)\, \Psi(t)\, +
  \int_0^t  \!\!  \phi(t-t') \, \lt[
 \int_{-\infty}^{+\infty}\!\!  w(x-x')\, p(x',t')\, dx'\rt]\,dt'\,,
\label{3.1BH}  %%\eqno(3.1)  $$
\ee
 in which 
$\delta (x)$ denotes the Dirac generalized function, and
 the {\it survival function}
 \bee %%$$
 \Psi(t) = \int_t^\infty \phi(t') \, dt' 
 \label{3.2BH}  %%\eqno(3.2)$$
 \ee
denotes the probability that at instant $t$ the particle
 is still sitting in its starting position
$x=0\,. $
%%%%%%%%%%%%%%
Clearly, Equation (\ref{3.1BH}) satisfies the initial condition
$p(x,0^+) = \delta (x)$.
%%%%%%%%%%%%
\vsp
Note that the {\it special choice}
\bee %% $$
w(x) = \delta (x-1)
\label{3.3BH} %% \eqno(3.3)$$
\ee
   gives the {\it pure renewal process}, with
position  $x(t)=N(t)$, denoting the {\it counting function},
and with jumps all of length 1 in positive
direction happening at the renewal instants.
\vsp
For many purposes the integral equation (\ref{3.1BH}) of CTRW
can be easily treated
by using the Laplace and Fourier transforms.
%%%%%%%%
Writing these as
$$     {\L} \lt\{ f(t);s\rt\}=  \widetilde f(s)
 := \int_0^{\infty} \!\!\e^{\ds \, -st}\, f(t)\, dt \,, \quad
 {\F} \lt\{g(x);\kappa\rt\}=  \widehat g(\kappa)
  := \int_{-\infty}^{+\infty} \!\! \e^{\,\ds +i\kappa x}\,g(x)\, dx
     \,,
$$
then in the Laplace-Fourier domain
 Equation (\ref{3.1BH})  reads
\bee %%$$ 
\widehat{\widetilde p}(\kappa ,s)
 =   {1-\widetilde\phi(s)  \over s} +
 \widetilde \phi(s) \,  \widehat w(\kappa )\,
   \widehat{\widetilde p}(\kappa ,s) \,.
   \label{3.4BH} %% \eqno(3.4)$$
   \ee
Introducing formally in the Laplace domain the auxiliary function 
\bee %% $$ 
\widetilde{H} (s) =
\frac{1- \widetilde{\phi}(s) }{ s\, \widetilde{\phi}(s)}
= \frac{\widetilde{\Psi}(s) }{\widetilde{\phi}(s)}
\,, \q \hbox{hence}\q
   \widetilde{\phi}(s) =  \frac{1}{1+s \widetilde{H} (s)} \,,
\label{3.5BH} %% \eqno(3.5)$$
\ee   
 and assuming that its Laplace inverse $H(t)$ exists, we get,
 following Mainardi et al. \cite{Mainardi PhysicaA00}, 
in the Laplace-Fourier domain the  equation 
\bee %% $$  
\widetilde{H} (s) \, \left[
  s\widehat{\widetilde p}(\kappa ,s)-1\right] =
  \left[ \widehat w(\kappa )-1\right]\,
   \widehat{\widetilde p}(\kappa ,s)
\,, 
\label{3.6BH} %%\eqno(3.6)$$
\ee
and in the space-time  domain  the generalized Kolmogorov-Feller equation
\bee %% $$ 
\int_0^t   H(t-t')\,
 \frac{\d}{\d t'} p(x,t')\, dt'  =
      -p(x,t) + \int_{-\infty}^{+\infty} w (x-x')\,p(x',t)\,dx',
\label{3.7BH} %%\eqno(3.7)$$
\ee
with $p(x,0) =\delta (x)$,
where $H(t)$ acts as a {\it memory function}.
\vsp
If  the Laplace inverse $H(t)$ of the formally introduced function $\widetilde H(s)$ does not exist, 
we can formally set
$\widetilde K(s) = 1/\widetilde H(s)$ and multiply 
 (\ref{3.6BH}) with $\widetilde K(s)$.
Then, if $K(t)$ exists, we get in place of (\ref{3.7BH}) the alternative
form of the generalized Kolmogorov-Feller equation
\bee %% $$ 
 \frac{\d}{\d t} p(x,t)\  =
     \int_0^t K(t-t')\,\left[ -p(x,t') + \int_{-\infty}^{+\infty} w (x-x')\,p(x',t')\,dx'
	 \right]\, dt'\,,
\label{3.7'BH} %% \eqno(3.7')$$
\ee
with $p(x,0) =\delta (x)$.
 where $K(t)$ acts as a {\it memory function}

%\vsp
Special choices of the memory function $H(t)$ are
$  {\mathbf{(i)}}$ and ${\mathbf{(ii)}}$, see Equations 
(\ref{3.8BH}) and (\ref{3.12BH}):
\vsp
\bee %%$$  
{\mathbf{(i)}}  \quad H(t)   =  \delta (t)
\q \hbox{corresponding to} \q \widetilde H(s)= 1\,, 
\label{3.8BH} %%\eqno(3.8)$$
\ee
%% and to the limiting value $\beta =1$,
 giving   the {\it exponential waiting time} with
\bee %% $$ 
\widetilde \phi (s) = {\ds \rec{1+s}}, \quad
 \phi (t) =    \Psi(t) = \e^{\ds \, -t}\,.
 \label{3.9BH} %%\eqno(3.9)  $$
\ee
In this case we obtain in the Fourier- Laplace domain
\bee %%$$ 
  s\widehat{\widetilde p}(\kappa ,s)-1 =
  \lt[ \widehat w(\kappa )-1\right]\,
   \widehat{\widetilde p}(\kappa ,s)\,, 
   \label{3.10BH}  %%\eqno(3.10)$$
   \ee
and in the space-time domain 
the {\it classical Kolmogorov-Feller equation}
\bee %% $$ 
\frac{\d}{ \d t} p(x,t) = - p(x,t)  +
   \int_{-\infty}^{+\infty}   w  (x-x')\,p(x',t)\,dx'\,,
     \quad p(x,0) = \delta(x)\,.
\label{3.11BH}%%\eqno(3.11) $$
\ee
\vsp
\bee %%$$ 
  {\mathbf{(ii)}}\quad H(t)   = {\ds \frac{t^{-\beta}}{ \Gamma(1-\beta)}}\,,\; 0<\beta <1\,,
\; \hbox{corresponding to} \; \widetilde H(s)= s^{\beta-1}\,,
\label{3.12BH}  %%\eqno(3.12)  $$
\ee
giving  the {\it Mittag-Leffler waiting time}
with
\bee %% $$ 
\widetilde \phi (s) = {\ds \rec{1+s^\beta }}, \quad
 \phi (t) =       - {\ds \frac{d}{dt} E_\beta (-t^\beta )}
 = \phi^{ML}(t),
 \quad
 \Psi(t) =E_\beta (-t^\beta)\,.
 \label{3.13BH} %%\eqno(3.13)  $$
 \ee
 %% then the time fractional continuous time random walk obeying the
In this case we obtain in the Fourier-Laplace domain
\bee %%$$ 
s^{\beta-1} \, \left[
  s\widehat{\widetilde p}(\kappa ,s)-1\rt] =
  \lt[ \widehat w(\kappa )-1\right]\,
   \widehat{\widetilde p}(\kappa ,s) \,, 
\label{3.14BH} %%   \eqno(3.14)$$
\ee
and in the space-time domain the {\it time fractional Kolmogorov-Feller equation} 
%% with the {\it Caputo} fractional time derivative
\bee %%$$ 
  \, _*D_t^\beta \,  p(x,t) =
     -  p(x,t) +   \int_{-\infty}^{+\infty} w(x-x')\,
   p(x',t) \, dx'\,, \quad p(x,0^+) = \delta(x)\,,  
   \label{3.15BH} %% \eqno(3.15)$$
   \ee
   where $\, _*D_t^\beta \,$ denotes the fractional derivative of of order $\beta$
   in the Caputo sense, see Section 3. %% \footnote{%%
%%%%%%%%%%%%%%%%%%%%%%%%%%%%%%
\vsp
The time fractional Kolmogorov-Feller equation can be also expressed via
the Riemann-Liouville fractional derivative $\,D_t^{1-\beta}$, see again
Section 3, 
that is
\bee %% $$ 
\frac{\d}{\dt}p(x,t) = \,D_t^{1-\beta}\,\left[
 -  p(x,t) +   \int_{-\infty}^{+\infty} w(x-x')\,
   p(x',t) \, dx'\right],   
   \label{3.16BH}  %%  \eqno(3.16)$$
   \ee
   with  $p(x,0^+) = \delta (x)$.
 The equivalence of the two forms (\ref{3.15BH}) and 
 (\ref{3.16BH}) is easily proved
 in the  Fourier-Laplace domain by multiplying both sides of 
 Equation (\ref{3.14BH})
 with the factor $s^{1-\beta}$.   
 %%  For this purpose we must rewrite  Equation (3.6)  in  
  %%  $$   \left[
  %% s\widehat{\widetilde p}(\kappa ,s)-1\rt] =
  %% \widetilde K(s)\, \lt[ \widehat w(\kappa )-1\right]\,
  %%  \widehat{\widetilde p}(\kappa ,s) \,,
  %%  \quad \widetilde K(s) = \frac{1}{\widetilde H(s)}\,, \eqno(3.15)$$
%%	and consequently assume the choice 
%%	$$   (iii) \quad K(t)   = {\ds \frac{t^{\beta-2}}{ \Gamma(\beta-1)}}\,,\; 0<\beta <1\,,
%% \q \hbox{corresponding to} \q \widetilde K(s)= {\ds s^{1-\beta}}
%% \eqno(3.16)  $$
%% In this case $K(t)$ is a generalized function so has no longer a meaning of a memory function.
\vsp
We  note that the choice ${\mathbf{(i)}}$  may be considered as a limit
of the choice ${\mathbf{(ii)}}$ as $\beta =1$.
In fact, in this limit we find   
$\widetilde H(s)  \equiv 1$ so 
$H(t)= t^{-1}/\Gamma(0)\equiv \delta(t)$
%(according to a formal representation of the Dirac generalized function
 %%\cite{GelfandShilov_64}),
so that Equations (\ref{3.6BH})--(\ref{3.7BH}) reduce to 
Equations (\ref{3.10BH})--(\ref{3.11BH}), respectively.
In this case the order of the Caputo  derivative
reduces to 1 and that of the R-L derivative to 0, 
whereas the Mittag-Leffler waiting time  law reduces to the
exponential.

%\vsp
In the sequel we will formally unite the choices 
({\bf i}) and ({\bf ii}) by defining what we call
the Mittag-Leffler memory function
\bee %%$$ 
 H^{ML}(t)= 
\left\{
\begin{array}{ll}
{\ds \frac{t^{-\beta}} {\Gamma(1-\beta)}} \,, &  \hbox{if} \q 0<\beta<1\,,
\\
{\ds \delta(t)}\,, & \hbox{if} \q \beta=1 \,,
\end{array}
\right .
\label{3.17BH} %%\eqno(3.17)$$
\ee
whose Laplace transform is
\bee %% $$
\widetilde H^{ML}(s)= s^{\beta -1}\,, \q 0<\beta \le 1\,.
\label{3.18BH} %%\eqno(3.18) $$
\ee
Thus we will consider the whole range $0<\beta \le 1$
 by  extending  the Mittag-Leffler waiting time law
in (\ref{3.13BH}) to include  the  exponential law 
(\ref{3.9BH}).  

\vsp
{\bf Remark:} Equation (\ref{3.7BH})
 clearly may be supplemented by an arbitrary
initial probability density $p(x,0)= f(x)$.
The corresponding  replacement of $\delta(x))$
by $f(x)$ in (\ref{3.1BH}) 
then requires in (\ref{3.4BH}) multiplication of the term
$(1-\widetilde\phi(s))/s$ by $\widehat f(\kappa)$ and 
in (\ref{3.6BH}) replacement of the
LHS by 
$\widetilde{H} (s) \, \left[
  s\widehat{\widetilde p}(\kappa ,s)-\widehat f(\kappa)\right]$. 
With $p(x,0) =\delta(x)$ we obtain in $p(x,t)$ the fundamental solution of 
Equation~(\ref{3.7BH}).
\vsp
{\bf Note:}
The probability density function for the waiting time distribution
 in terms of the Mittag-Leffler function was formerly given  since 1995 by Hilfer
\cite{Hilfer FRACTALS95}, \cite{Hilfer PhysA03}, \cite{Hilfer-Anton PRE95}.
In these papers  the waiting time density  was given with the Mittag-Leffler function in two parameters without noting    the relation with the first derivative of the classical Mittag-Leffler function as stated in Equation (\ref{1.45}). 
We also note that 10 years earlier
Balakrishnan \cite{Balakrishnan PhysA85} 
had derived a similar expression without
recognizing the Mittag-Leffler function. Like in the case of the thinning process 
dealt by Gnedenko-Kowalenko (see section 7) once again the Mitag-Leffler function
was unknown for the authors.

\subsection{Manipulations: Rescaling and Respeeding}  %% SECTION 4

We now consider two types of manipulations on the CTRW
by acting  on its  governing Equation~(\ref{3.7BH}) 
in its Laplace-Fourier representation (\ref{3.6BH}).
\\
({\bf A}): rescaling the waiting time, hence the whole time axis;
\\
({\bf B}): respeeding the process.
\vsp
({\bf A}) means change of the unit of time (measurement).
We replace the random waiting
time  $T$  by a waiting time $\tau T$,
with the positive {\it rescaling factor} $\tau$.
Our idea is to take $0<\tau \ll 1$   in order to bring into near
sight the distant future.
In a moderate span of time we will so have a large number of jump
events. For $\tau >0$  we get the rescaled waiting time density
\bee %% $$\phi _\tau (t) = \phi (t/\tau )/\tau\,, \q \hbox{hence} \q
\widetilde\phi _\tau (s) = \widetilde\phi (\tau s)\,.  
\label{4.1BH} %%\eqno(4.1)$$
\ee
By decorating also the density $p$ with an index $\tau $   we
obtain the rescaled integral equation of the CTRW
in the Laplace-Fourier  domain as
\bee %% $$ 
\widetilde{H}_\tau  (s) \, \lt[
  s\widehat{\widetilde {p}}_\tau(\kappa ,s)-1\rt] =
  \lt[ \widehat w(\kappa )-1\rt]\,
   \widehat{\widetilde {p}}_\tau (\kappa ,s)\,,
\label{4.2BH} %%\eqno(4.2)$$
\ee
where, in analogy to (\ref{3.5BH}),
\bee %% $$ 
\widetilde{H}_\tau (s) =
\frac{1- \widetilde{\phi}(\tau s) }{ s\, \widetilde{\phi}(\tau s)}\,.
\label{4.3BH} %%\eqno(4.3)$$
\ee
\vsp
({\bf B})  means multiplying the quantity representing ${\ds \frac{\d}{\dt} p (x,t)}$
by a factor $1/a$, where $a>0$ is the {\it respeeding factor}:
%% \\
$a>1$ means {\it acceleration}, $0<a<1$  means {\it deceleration}.
In the Laplace-Fourier representation 
this means multiplying the RHS of Equation (\ref{3.6BH}) by the factor $a$  
since the expression
 $\lt[s\widehat{\widetilde {p}}(\kappa ,s)-1\rt]$
    corresponds to
$ {\ds \frac{\d}{\dt} p (x,t)}$.
%% Note  that $\widehat{p}_\tau(\kappa ,t=0)= \widehat{\delta}(\kappa)=1\,$
\vsp
We now chose
to consider the procedures of rescaling and respeeding
in their combination so that the equation in the transformed domain 
of the  rescaled and respeeded process has the form
\bee %% $$  
\widetilde{H}_\tau  (s) \, \lt[
  s\widehat{\widetilde {p}}_{\tau,a}(\kappa ,s)-1\rt] =
 a\, \lt[ \widehat w(\kappa )-1\rt]\,
   \widehat{\widetilde {p}}_{\tau,a} (\kappa ,s)\,, 
\label{4.4BH}  %% \eqno(4.4)$$
\ee
%% \vsp
Clearly, the two manipulations  can be discussed separately:
the choice $\{\tau >0,\, a=1\}$ means {\it pure rescaling},
the choice $\{\tau =1,\, a >0\}$ means {\it pure respeeding}
of the original process.
In the special case  $\tau =1$   we only respeed the
original system; if  $0 <\tau \ll 1$
  we can counteract the compression effected by rescaling to again obtain a
moderate number of events in a moderate span of time by respeeding
(decelerating) with  $0<a \ll 1$.
 These vague notions will become clear as soon as we consider power law
waiting times.
\vsp
Defining now
\bee % $$ 
\widetilde{H}_{\tau,a}  (s) := \frac{\widetilde{H}_{\tau}  (s)}{a}=
\frac{1- \widetilde{\phi}(\tau s) }{ as\, \widetilde{\phi}(\tau s)}\,.
\label{4.5BH} %%\eqno(4.5)$$
\ee
we finally get, in analogy to (\ref{3.6BH}), the equation
\bee %% $$ 
\widetilde{H}_{\tau,a}  (s) \, \lt[
  s\widehat{\widetilde {p}}_{\tau,a}(\kappa ,s)-1\rt] =
  \lt[ \widehat w(\kappa )-1\rt]\,
   \widehat{\widetilde {p}}_{\tau,a}(\kappa,s)\,.
\label{4.6BH}  %%   \eqno(4.6)$$
\ee
%%%%%%%%%%%%%%%%%%
What is the combined effect of rescaling and respeeding on the waiting
time density?
\vsp
In analogy to (\ref{3.5BH}) and taking account of (\ref{4.5BH})  we find
\bee %% $$
   \widetilde{\phi}_{\tau ,a}(s) =
        \frac{1}{1+s \widetilde{H}_{\tau ,a} (s)}=
        \frac{1}
{1+s {\ds\frac{1-\widetilde{\phi}(\tau s)}{as \,\widetilde{\phi}(\tau s)}}}
\,,
\label{4.7BH} %%\eqno(4.7)$$
\ee
and so, for the deformation of the waiting time density, the
{\it essential formula}
\bee %% $$   \widetilde{\phi}_{\tau ,a}(s) =
     \frac{a\,\widetilde{\phi}(\tau s) }
     {1- (1-a)\widetilde{\phi}(\tau s)}
%%   = \frac{a \, \widetilde{\phi}(\tau s)}
%%     {a\widetilde{\phi}(\tau s)+1 -\widetilde{\phi}(\tau s)}
\,. 
\label{4.8BH} %%\eqno(4.8) $$
\ee
\vsp
{\bf Remark}:
The formula (\ref{4.8BH}) has the same structure as the thinning
formula (\ref{2.5BH})  in Section 5 (just devoted to the thinning theory) 
by  identification of  $a$  with $q$. 
In both problems  we have a rescaled process defined by a time scale $\tau$,
and  we send the  relevant factors $\tau$, $a$ and $q$ to zero
under a proper relationship.
However in the thinning theory the relevant independent parameter going to 0 
is that of thinning (actually respeeding) whereas in the present problem 
it is the rescaling parameter $\tau$.

%%%%%%%%%%%%%%%%%%%%%%%%%%%%%%%%%%%%%

%\section{the Universality of the Mittag-Leffler Function as Probability Density Function}

\section{Power Laws and Asymptotic Universality of the Mittag-Leffler Waiting
Time Density}    %%% SECTION 5 BH

We have essentially two different situations for the waiting time distribution
according to its first moment (the expectation value)  being  finite or infinite.
In other words we assume for the waiting time $pdf$ $\phi(t)$
 either
 \bee %% $$ 
\rho  := \int_0^\infty t'\,  \phi(t')\, dt' < \infty \,, \quad
\hbox{labelled as}\; \beta=1\,,  
\label{5.1BH} %%\eqno (5.1)$$
\ee
%% In the second case (dispensing with decorating by a slowly varying function at infinity)
%% we have the following asymptotic power law behaviours
%% as $t \to \infty$: 
or
\bee %% $$  
\phi(t) \sim c \, t^{-(\beta +1)} \; \hbox{for}\; t \to \infty \quad \hbox{hence} \;
 \; \Psi(t) \sim \frac{c}{\beta} \, t^{-\beta } \,,
 \; 0<\beta<1\,,\; c>0\,. 
 \label{5.2BH}  %%\eqno(5.2)$$
 \ee
 For convenience we have dispensed in 
 (\ref{5.2BH}) with decorating by a slowly varying function at 
 infinity  the asymptotic power law.
Then, by the standard Tauberian theory (see References \cite{Feller BOOK71,Widder BOOK46})
the above  conditions (\ref{5.1BH})--(\ref{5.2BH}) mean in the Laplace domain
 the (comprehensive) asymptotic  
form
\bee %% $$  
\widetilde{\phi}(s) =1 -\lambda s^\beta  + o(s^\beta )
\q \hbox{for} \q s \to 0^+\,,
\q  %%% \hbox{with}  \q 
0<\beta \le 1\,,
\label{5.3BH} %% \eqno(5.3)$$
\ee
where we have 
\bee %% $$
  \lambda = \rho\,,\q \hbox{if}\q \beta=1\,; \; 
 \lambda = c \Gamma(-\beta )=
\frac{c}{\Gamma (\beta+1)} \, \frac{\pi} {\sin(\beta \pi)}\,,
\; \hbox{if}\; 0<\beta<1\,.
 \label{5.4BH} %% \eqno(5.4) $$
 \ee
Then, {\it fixing $s$}    as required by the continuity theorem of
probability theory
for Laplace transforms, taking
\bee %% $$
a =\lambda \tau ^\beta\,,
\label{5.5BH}   %%\eqno(5.5)$$
\ee
 and {\it sending $\tau $  to zero}, we obtain in the limit
the Mittag-Leffler waiting time law.
In fact,  Equations (\ref{4.8BH}) and (\ref{5.3BH}) imply  as $\tau \to 0$
with  $0<\beta \le 1$,
\bee %%$$    
\widetilde{\phi}_{\tau ,\lambda \tau ^\beta} (s) 
=  \frac{\lambda \tau^\beta\, \left[1-\lambda \tau^\beta s^\beta + o(\tau^\beta s^\beta)\right]}
{1 - (1-\lambda \tau^\beta)\, \left[1-\lambda \tau^\beta s^\beta + o(\tau^\beta s^\beta)\right] }
  \to
\frac{1}{1 + s^\beta }\,, 
\label{5.6BH} %% \eqno(5.6)$$
\ee
the Laplace transform of $\phi^{ML}(t)$. %% see (1.1).%%  and Appendix C.
This formula expresses {\bf the asymptotic universality of
the Mittag-Leffler waiting time law} that includes the exponential law for $\beta=1$.
It can easily be generalized to the case of power laws decorated with
slowly varying functions, thereby using the Tauberian theory by Karamata 
(see again References \cite{Feller BOOK71,Widder BOOK46}).
\vsp
{\bf Comment:} The formula (\ref{5.6BH}) says that our general power law waiting time
density is
gradually deformed into the Mittag-Leffler waiting time density as
$\tau $ tends to zero.
\vsp
{\bf Remark:} Let us stress here the distinguished character of the Mittag-Leffler
waiting time density ${\ds \phi^{ML}(t)= - \frac{d}{dt} E_\beta(-t^\beta)}$.
%%defined in (1.1).
Considering its Laplace transform   
\bee %% $$  
 \widetilde{\phi}^{ML}(s) = \frac{1}{1+s^\beta }\,, \q
   \phi^{ML} (t) =  - {\ds \frac{d}{dt} E_\beta (-t^\beta )}\,,
 \; 0<\beta \le 1\,,
\label{5.7BH} %%\eqno(5.7)$$
\ee
we can easily prove the identity
\bee %% $$
    \widetilde{\phi}^{ML}_{\tau,a} (s)
%% :=\frac{\widetilde{\phi}^{ML (\tau s)}{a}
= \widetilde{\phi}^{ML} (\tau s/a^{1/\beta })
\q \hbox{for all} \q \tau >0, \q a>0\,. 
\label{5.8BH}%%\eqno(5.8)$$
\ee
Note that Equation (\ref{5.8BH}) states the 
{\it self-similarity} of  the combined operation
{\it rescaling-respeeding}
for the Mittag-Leffler waiting time density. 
In fact, (\ref{5.8BH}) implies
$ {\phi}^{ML}_{\tau,a} (t) =
   {\phi}^{ML}(t/c) /c$ with
$ c= \tau /a^{1/\beta }\,,$
 which means replacing the random waiting time $T^{ML}$
by   $ c\, T^{ML}$.
As a consequences, choosing  $a= \tau ^\beta $ we have
\bee %% $$    
\widetilde{\phi}^{ML}_{\tau,\tau^\beta} (s)
%% :=\frac{\widetilde{\phi}^{ML (\tau s)}{\tau^\beta}
 = \widetilde{\phi}^{ML} ( s)
\q \hbox{for all} \q \tau >0\,. 
\label{5.9BH} %%\eqno(5.9)$$
\ee
Hence {\it the Mittag-Leffler waiting time density  is invariant against
    combined rescaling  with $\tau$ and respeeding with  
    $a=\tau^\beta$}.
	
%	\vsp
	Observing  (\ref{5.6BH}) we can say that $\phi^{ML} (t)$ is a $\tau \to 0$
	attractor for any power law waiting time
	(\ref{5.2BH}) under simultaneous rescaling with $\tau$ and respeeding with
	$a= \lambda \tau^\beta$.
%% \vsp {\bf Remark}:	
In other words, this attraction property of the Mittag-Leffler probability 
 distribution with respect to power law  waiting  times (with $0<\beta\le  1$) 
 is a kind of analogy to the attraction of sums of power law jump 
 distributions by stable distributions.
%%%%%%
\section{The Mittag-Leffler functions w.r.t. the time fractional diffusion-wave equations and the Wright functions}
In this section we  show the relations of the Mittag-Leffler function
with the Wright function via Laplace and Fourier transformations, in order to provide 
other arguments to outline  the role of the Mittag-Leffler in the Fractional Calculus.
For this purpose, because of the necessity  to work with  two independent parameters 
 we first recall
 the proper definitions of the Mittag-Leffler and the Wright function.
 Then we will consider the time fractional diffusion-wave equation with its
 fundamental solutions to the basic boundary value problem that urn out to be expressed in terms of some special cases of the Wright functions, the so called $F$ and $M$ functions.
 Finally we pay attention to some noteworthy formulas for the $M$-Wright function, including 
 its connections with the Mittag-Leffler function.
 \vsp
 \subsection{Definitions and Main Properties  of the Wright Functions}

%\vsp
The classical \emph{Wright function}, that we denote by $W_{\lambda , \mu}(z)$, is defined by the series representation convergent in the whole complex plane,
\bee
W_{\lambda , \mu}(z) :=\sum_{n=0}^{\infty}{\frac{z^n}{n!\Gamma (\lambda n + \mu)}}, ~~~ \lambda > -1, ~~~ \mu \in \mathbb{C},
\label{F.1}
\ee
As a consequence 
$W_{\lambda , \mu}(z)$ is an \emph{entire function} for all
$\lambda \in (-1, +\infty)$. 
 Originally Wright assumed $\lambda \ge 0$ in connection with his investigations on the asymptotic theory of partition
 \cite{Wright 33,Wright 35} and only in 1940 he considered
  $-1 < \lambda < 0$, \cite{Wright 40}.
We note that in the Vol 3, Chapter 18 of the handbook} of the Bateman Project
\cite{Erdelyi BATEMAN}, presumably for a misprint,
 the parameter $\lambda$  is restricted to be non-negative,
whereas  the Wright functions remained  practically ignored   in other  handbooks.  
In 1993  the present author,   being aware only  of  the Bateman handbook,
proved that the Wright function is entire also for $-1<\lambda<0$  in his approaches to the time fractional diffusion equation, as outlined in his papers
published  from 1994 to 1997,
\cite{Mainardi WASCOM93,Mainardi RADIOPHYSICS95,Mainardi AML96,%%
Mainardi CHAOS96,Mainardi CISM97}.
For other earlier  treatments of this function we refer to the 
1999 paper by Gorenflo, Luchko and Mainardi \cite{GOLUMA 99}).
%%%

%\vsp
In view of the asymptotic representation in the complex domain 
and of the Laplace transform 
%for positive argument  $z=r>0$
  the Wright functions were  distinguished by the author in
   \emph{first kind} ($\lambda \geq 0$) and \emph{second kind}
 ($-1< \lambda < 0$) 
 as outlined \eg in the Appendix F of his book 
 \cite{Mainardi BOOK10}.
 %% and more recently
 % in the survey by Mainardi and Consiglio
 %In particular,  for the asymptotic behaviour, we refer the interested reader to
% the two  papers by Wong and Zhao \cite{Wong 99a,Wong 99b},
 % and to  the surveys by Luchko and by Paris in the Handbook of Fractional Calculus
%  and Applications,
%   see respectively \cite{Luchko HFCA}, \cite{Paris HFCA},
%    and references therein.
 
% \vsp 
We note that the Wright functions  are  entire of order $1/(1+\lambda)$
 hence only  the first kind  functions ($\lambda \ge 0$)  are of exponential order whereas  the second kind functions ($-1<\lambda<0$)  
 are not of exponential\ order.
 The case $\lambda =0$ is trivial since
 $W_{0,\mu}(z) = {\e^z}/{\Gamma(\mu)}.$

%\vsp
Following the profs in the Appendix F in Reference \cite{Mainardi BOOK10} we get
 the following Laplace transform pairs of the Wright functions in terms of
 the Mittag-Leffler functions in two parameters, where
 $r$ can be the time variable $t>0$ or the space variable $x>0$)
 \\
 {\it for the first kind} ($\lambda \ge 0$)
 \bee %% $$  
  W_{\lambda,\mu } (\pm r) \,\div \,
    \rec{s}\, E_{\lambda ,\mu }\left(\pm \rec{s}\right) \,,
  \q \lambda > 0\,, %% \q |s| > \rho> 0\,, 
  \label{F.22} %%\eqno(F.22)$$
  \ee
  \\
  {\it for the second kind} ($\lambda = -\nu, \; 0<\nu<1$)
  \bee  %%  $$
     W_{-\nu ,\mu }(-r) \,\div\,
   E_{\nu  , \mu +\nu  }(-s)\,, \q 0<\nu  <1\,.
   \label{F.25} %% \eqno(F.25)$$
  \ee    
The Wright functions of the first kind are useful to find the solutions 
of some (linear and non-linear) differential equations of fractional order as recently 
shown by Garra and Mainardi, \cite{Garra-Mainardi ROMP20}.

%\vsp
Since the pioneering works in 1990's by the author,  noteworthy cases of Wright functions 
of the second kind, known as {\it auxiliary functions} $F$ and $M$ play   
fundamental roles in solving  the  Signalling problem
and the Cauchy value problem, respectively
for the time fractional diffusion-wave equation.

%\vsp
 We first recall hereafter these auxiliary functions in terms of the Wright functions of the second kind, following their power series representations. They read
\bee %%$$ 
F_\nu (z) :=   W _{-\nu , 0}(-z)\,, \q 0<\nu<1\,, 
\label{F.9} %%\eqno(F.9)$$
\ee
and
\bee  %% $$
 M_\nu (z) :=  W _{-\nu , 1-\nu }(-z)\,,
\q 0<\nu<1 \,,
\label{F.10}  %%  \eqno(F.10)$$
\ee
interrelated through
\bee %% $$ 
F_\nu (z) = \nu  \, z \, M_\nu (z ) \,.
\label{F.11}  %%\eqno(F.11)$$
\ee
%% As it is shown in Chapter 6, the motivation was  based on the inversion of certain Laplace %transforms  in order to obtain the  fundamental solutions of the fractional diffusion-wave %equation in the space-time domain. 
\vsp
The {\it series representations}  of our auxiliary functions are derived 
from those of $W_{\lambda, \mu}(z)$ in (\ref{F.1}). We have:
\bee %% $$ 
 %%\begin{array}{ll}
 F_\nu (z) =
 {\ds \sum_{n=1}^{\infty}
\frac{(-z)^n}{  n!\, \Gamma(-\nu n)}}  
= -{\ds \rec{\pi}\, \sum_{n=1}^{\infty}
 \frac{(-z)^{n}}{ n!}\,
 \Gamma(\nu n +1 )\, \sin(\pi \nu  n)\,,}
%\end{array}
\label{F.12} %%\eqno(F.12)   $$
\ee
 and
\bee  %%$$ 
%\begin{array}{ll}
M_\nu (z) =
 {\ds \sum_{n=0}^{\infty}
 \frac{(-z)^n }{  n!\, \Gamma[-\nu n + (1-\nu )]} }  
   = {\ds \rec{\pi}\, \sum_{n=1}^{\infty}\,\frac{(-z)^{n-1} }{  (n-1)!}\,
  \Gamma(\nu n)  \,\sin (\pi\nu n)}  \,,
%% \end{array}
\label{F.13} %%   \eqno(F.13)$$
\ee  
where we have used the well-known reflection formula for the Gamma function,
 $$ \Gamma(\zeta)\,\Gamma(1-\zeta)  =\pi /\sin\,\pi \zeta\,.$$
%% \vsp
%%
\subsection{The Time-Fractional Diffusion-Wave Equation and the Related Green Functions}
For the reader's convenience let us recall the main formulas for the time fractional diffusion equations  and their fundamental solutions (also referred to as the  Green functions) for the Cauchy and Signalling problems.
For more details we refer  to References \cite{Mainardi LUMAPA01},
\cite{Luchko-Mainardi HFCA19}.

%\vsp
 Denoting as usual   $x,t $ the space and time variables, 
 and $r=r(x,t)$ the response variable, 
  the family of these evolution equations reads  
\bee %%  $$
{\d^{\beta} r\over \dt^{\beta} } =
  a\,{\d^2 r \over \dx^2}  \,,	\q 0<\beta \le 2\,,
\label{6.3} %%\eqno(6.3)$$
\ee
where {\it the time derivative of order $\beta$ is intended in the Caputo sense},
namely is the  operator $\,_*D_t^\beta$, introduced in Section 3, but 
for order less than 1, see Equation (\ref{2.6b}),
and $a$ is a positive
constant of dimension $L^2 \, T^{-\beta}$.
%%%%%
Thus  we must distinguish
the cases $0<\beta\le 1$ and $1<\beta\le 2$. We have
\bee %%$$ 
 {\d^{\beta} r\over \dt^{\beta} } :=  
\left\{
%%%\begin{cases}
\begin{array}{ll}
  {\ds \rec{\Gamma(1-\beta )}} \,
  {\ds \int_0^t} \!\!  \lt[{\ds{\d \over \d \tau}} \, r(x,\tau)\rt]\,
  {\ds{d\tau \over (t-\tau)^{\beta}}}  \,, & \, 0<\beta  < 1\,,\\ \\
{\ds {\d r\over \d t}}\,, & \; \beta =1\,; \\
%\end{cases}
\end{array}
\right .
\label{6.4a} %%\eqno (6.4a)$$ 
\ee
%%%%%%
\bee %%$$
 {\d^{\beta} r\over \dt^{\beta}} \!:=\! 
\left\{
 %\begin{cases}
 \begin{array}{ll}
 {\ds\rec{\Gamma(2-\beta )}} 
  {\ds \int_0^t}  \!\! \lt[{\ds{\d^2 \over \d \tau^2}} \, r(x,\tau)\rt]\,
  {\ds{d\tau \over (t-\tau)^{\beta -1}}} ,
 &  1<\beta < 2,\\ \\
 {\ds {\d^2 r\over \d t^2}}\,,  &  \beta =2\,.
%% \end{cases}
 \end{array}
\right . 
\label{6.4b} %% \eqno(6.4b)$$
\ee
%% where $\Gamma$ denotes the Gamma function.
\vsp
It should be noted that in both cases 
 $0<\beta\le 1, \;1<\beta \le 2$,
the  time fractional derivative in the L.H.S.
of  Equation (\ref{6.3})
can be removed by a suitable fractional integration.
leading to alternative forms where  the necessary initial   conditions at $t=0^+$
 explicitly appear.
 \vsp 
For this purpose we apply to Eq. (\ref{6.3}) the fractional integral 
operator of order $\beta$, namely 
$$
J_t^\beta f(t):= \frac{1}{\Gamma(\beta)}\, \int_0^t (t-\tau)^{\beta-1}\, f(\tau)\, d\tau.
$$
For  $\beta \in (0,1] $ we have:
 $$ J_t^\beta \,\circ \, _*D_t^{\beta}\, r(x,t)=
 \, J_t^\beta\,\circ \, J_t^{1-\beta}\, D_t^1\, r(x,t) =
  \,J_t^1\, D_t^1\, r(x,t) = r(x,t) - r(x,0^+)\,.$$
  For $\beta \in (1,2] $ we have:
  $$ J_t^\beta \,\circ \, _*D_t^{\beta}\, r(x,t)\!=\!
 \, J_t^\beta\,\circ \, I_t^{2-\beta}\, D_t^2\, r(x,t)\! =
\!  \,J_t^2\, D_t^2\, r(x,t) \!=\! r(x,t) -  r(x,0^+) -t\, r_t(x,0^+).)$$  
\vsp
Then, as a matter fact,  we get  the integro-differential equations:
\vsn if $0<\beta \le 1\,:$
\bee %%$$ 
 r(x,t) =
r(x,0^+)
 + {\ds{a\over \Gamma(\beta )}}\,
  {\ds \int_0^t} \!\! 
 \left({\ds{\d^2 r\over \dx^2}}\right) \,  (t-\tau)^{\beta -1}\, d\tau\,;
\label{6.5a}  %%   \eqno(6.5a)$$
\ee
%%%%%%
\vsn if $1<\beta \le 2\,:$
\bee %% $$ \!\! r(x,t) \!=\!
 r(x,0^+) +  t \,\frac{\d}{\d t}\left. r(x,t)\right|_{t=0^+}
 + {\ds{a\over \Gamma(\beta )}}
  {\ds\int_0^t} \!\!
 \left({\ds{\d^2 r\over \dx^2}}\right)  (t-\tau)^{\beta -1}\, d\tau.
\label{6.5b} %%    \eqno(6.5b)$$
\ee
%%%%%%%%%
Denoting by
$f(x)\,,\, x\in \RR $ and $h(t)\,,\, t\in \RR^+\,$
sufficiently
well-behaved functions,    the basic boundary-value problems are thus formulated as
following, assuming $0<\beta \le 1$,
\vsn
{\it (a) Cauchy problem}
\bee %% $$
   r(x,0^+)=f(x) \,,	\; -\infty <x < +\infty\,;
  \;	 r(\mp \infty,t) = 0\,,\; \, t>0\,;  
  \label{6.6a}  %%\eqno(6.6a)$$
  \ee
%%%%%%%%%%
\vsn
{\it (b) Signalling problem}
\bee %% $$ 
 r(x, 0^+)  =  0 \,, \;  x>0\,;\;
    r(0^+,t ) =h(t) \,,\; r(+\infty,t) =0 \,, \;   t >0
\,. 
\label{6.6b} %%\eqno(6.6b) $$
\ee
If $1 <\beta < 2\,, $ we must add into (\ref{6.6a}) and (\ref{6.6b})
the initial values of the first time derivative of the field variable,
$r_t(x,0^+)\,,$ since in this case
the  corresponding fractional derivative is expressed in terms
of the second order  time derivative.
%% MISLEADING according to Gorenflo
%% and, consequently, two linearly independent solutions are to be
%% determined.	However,
To ensure the continuous dependence  of our solution
with respect to the parameter $\beta  $
also in the transition from $\beta  =1^-$ to  $\beta =1^+\,,$
we agree  to  assume 
%% $$r_t(x,0^+)
\bee %% $$
\frac{\d}{\d t}\left. r(x,t)\right|_{t=0^+} = 0\,, \; 
\hbox{for} \; 1<\beta \le 2\,,
\label{6.7} %%\eqno (6.7)$$
\ee
as it turns out from the integral forms (\ref{6.5a})--(\ref{6.5b}).
%%%%%
\vsp
In view of our subsequent analysis we find it convenient to set
\bee %% $$ 
 \nu  :={\beta / 2}\,,
   \q {\hbox{so}}\q
   \begin{cases}
    0<\nu  \le 1/2 \,, \, \Longleftrightarrow \, 0<\beta\le 1\,,\\
	1/2 <\nu \le 1\,, \, \Longleftrightarrow \,  1<\beta \le 2\,,
	\end{cases} 
	\label{6.8} %% \eqno(6.8) $$
	\ee
and from now on to add the parameter  $\nu $ to the
independent space-time variables $x\,,\,t$ in the solutions,
writing $r = r(x,t;\nu )$.

%\vsn
For the {Cauchy} and {Signalling}  problems we
introduce
the so-called {\it Green functions} $\Gc (x,t;\nu )$ and $\Gs(x,t;\nu )$,
which represent the respective fundamental solutions,
obtained when $f(x) = \delta (x)$ and $h(t) = \delta (t)\,. $
As a consequence, the solutions of the two basic problems
are obtained by a space or time convolution  according to
\bee %% $$ 
r(x,t;\nu)
= \int_{-\infty}^{+\infty} \Gc(x-\xi ,t;\nu  ) \, f(\xi)
\, d\xi     \,,  
\label{6.9a}  %% \eqno(6.9a)$$
\ee
%%%%%%
\bee %% $$ 
r(x,t;\nu ) = \int_{0^-}^{t^+} \Gs(x,t-\tau;\nu  ) \, h(\tau ) \,
d\tau	  \,.  
\label{6.9b} %%\eqno(6.9b)$$
\ee
It should be noted that in (\ref{6.9a}) $\Gc(x ,t ;\nu) =
  \Gc(|x|,t ;\nu)$ because the Green function of the Cauchy problem turns
out to be an even function of $x$.
According  to a usual convention,
in (\ref{6.9b}) the limits of integration are extended to take into account
for the possibility of impulse functions centred at the extremes.
%%%%%%%

%\\
Now we recall the results obtained in 1990's by the author that allow us to express the two Green functions in terms of the auxiliary functions $F_\nu(\xi)$ and $M_\nu(\xi)$
where, for $x>0$,  $t>0$ 
\bee %% $$   
\xi :={x/(\sqrt{a}\, t^{\nu  })} >0
\label{6.22}  %% \eqno(6.22)$$
\ee
acts as
{\it similarity variable}.
Then   
we obtain the Green functions in the space-time domain in the form 
\bee %% $$ 
\Gc(x,t;\nu) = {1\over 2\, \nu\, x}\, F_\nu(\xi)
   =	{1\over2 \sqrt{a}\, t^\nu}\, M_\nu(\xi)
 \,,  
 \label{6.23a} %% \eqno(6.23a)$$
 \ee
 %%%%
 \bee %%$$ 
 \Gs(x,t;\nu) = {1\over t}\, F_\nu(\xi)
 =  {\nu \, x\over \sqrt{a} \,t^{1+\nu}}\, M_\nu(\xi)\,.
 \label{6.23b} %% \eqno(6.23b)$$
\ee
  We also recognize the following {\it reciprocity relation} 
  for the  original Green functions, 
\bee %% $$ 
2\nu  \, x\,  \Gc(x,t;\nu  )  = t\, \Gs(x,t;\nu  )
 =   F_\nu(\xi) = \nu  \xi\, M_\nu(\xi)
  \,.
 \label{6.24} %%  \eqno(6.24) $$ 
 \ee
    Now $F_\nu(\xi)$, $M_\nu(\xi)$
are the {\it auxiliary functions} for the general case $0<\nu\le 1$,
which generalize those well known for the standard (Fourier) diffusion equation 
and for the standard (D'alembert) wave equation 
derived for  $\nu=1/2$ and for $\nu=1$, respectively.
%%%%
\subsection{Some Noteworthy Results for the $M_ \nu$ Wright Function}
In this survey we find  worthwhile to concentrate our attention 
on a single auxiliary function, the $M$-Wright function, 
sometimes referred to as the {\it Mainardi function}.
Indeed this function is indeed referred 
with this name 
in n the 1999 book by Podlubny \cite{Podlubny BOOK99}, that is one of the most cited treatises on fractional calculus.
 Then   this name is found in several successive  papers and books related
 to fractional diffusion and wave processes, 
 see for example, the relevant 2015 paper by Sandev et al. \cite{Sandev PRE15}.

Let us now  recall some interesting analytic results related to the so-called
Mainardi function.
 %%%
  One reason for the major attention is due to its straightforward generalization of the Gaussian probability density obtained for $\nu=1/2$, that is the fundamental solution of the Cauchy problem  for the standard diffusion equation.
 Furthermore it allows an impressive  visualization of the evolution with the order 
 $\nu \in (0,1)$ of 
 the Green function of the Cauchy problem of the fractional diffusion wave Equation
 (\ref{6.23a}) as shown 
in  Fig. 6 and Fig. 7 
  with $a=1$ and taking $t=1$.

 \begin{figure}   [h!]
	\centering
	\includegraphics[width= 100mm, height=60mm]{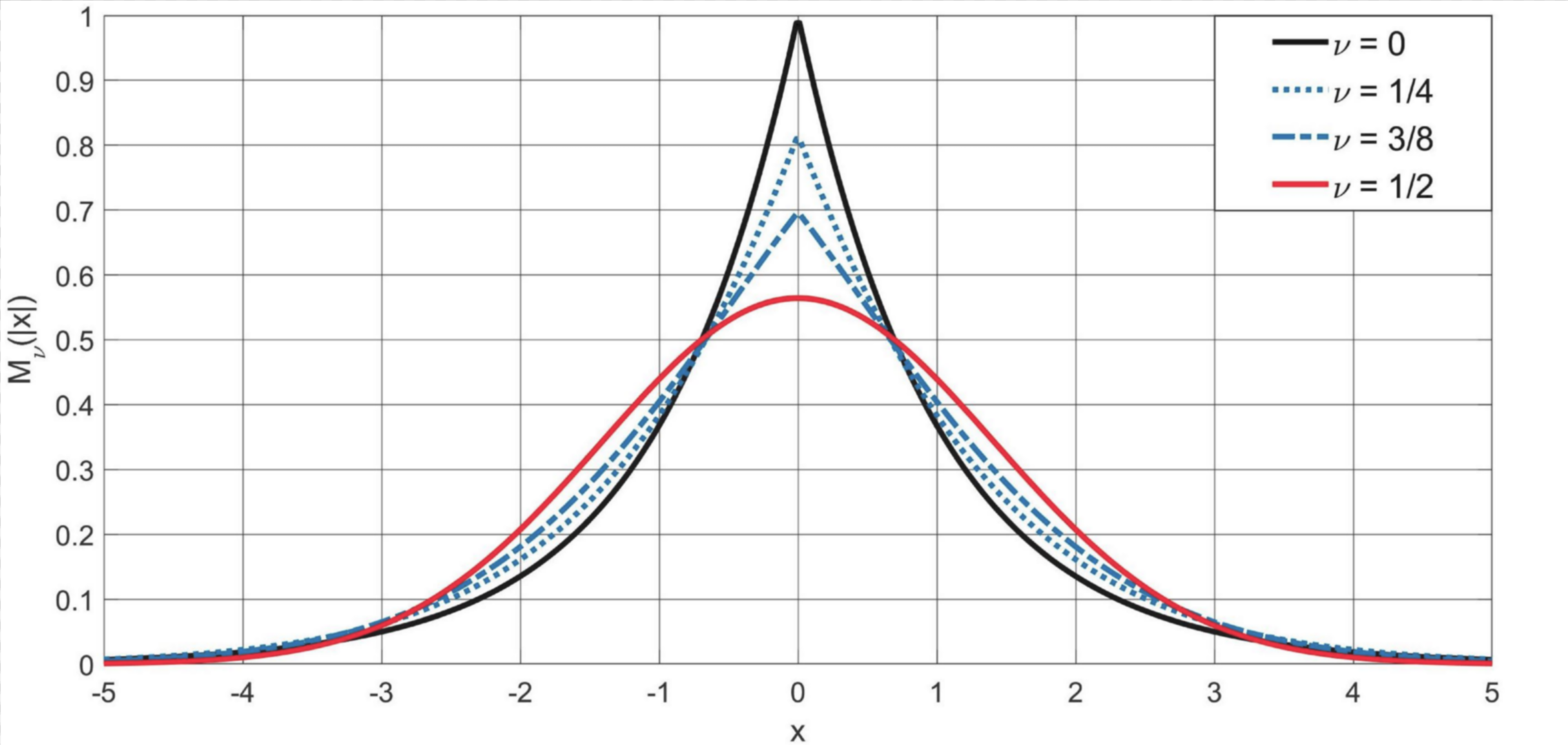}
%	\caption 
\vskip -0.45truecm
\caption{Plot of the symmetric $M-$Wright  function $M_{\nu}(|x|)$ 
	for $0 \leq \nu \leq 1/2$. Note that the $M-$Wright function becomes
a Gaussian density  for $\nu=1/2$.}	
\label{fig6}
	 \end{figure}%please cite figure 6 and figure 7 in the main text.
	 %%%%
\begin{figure}     [h!]
	\centering
	\includegraphics[width= 100mm, height=60mm]{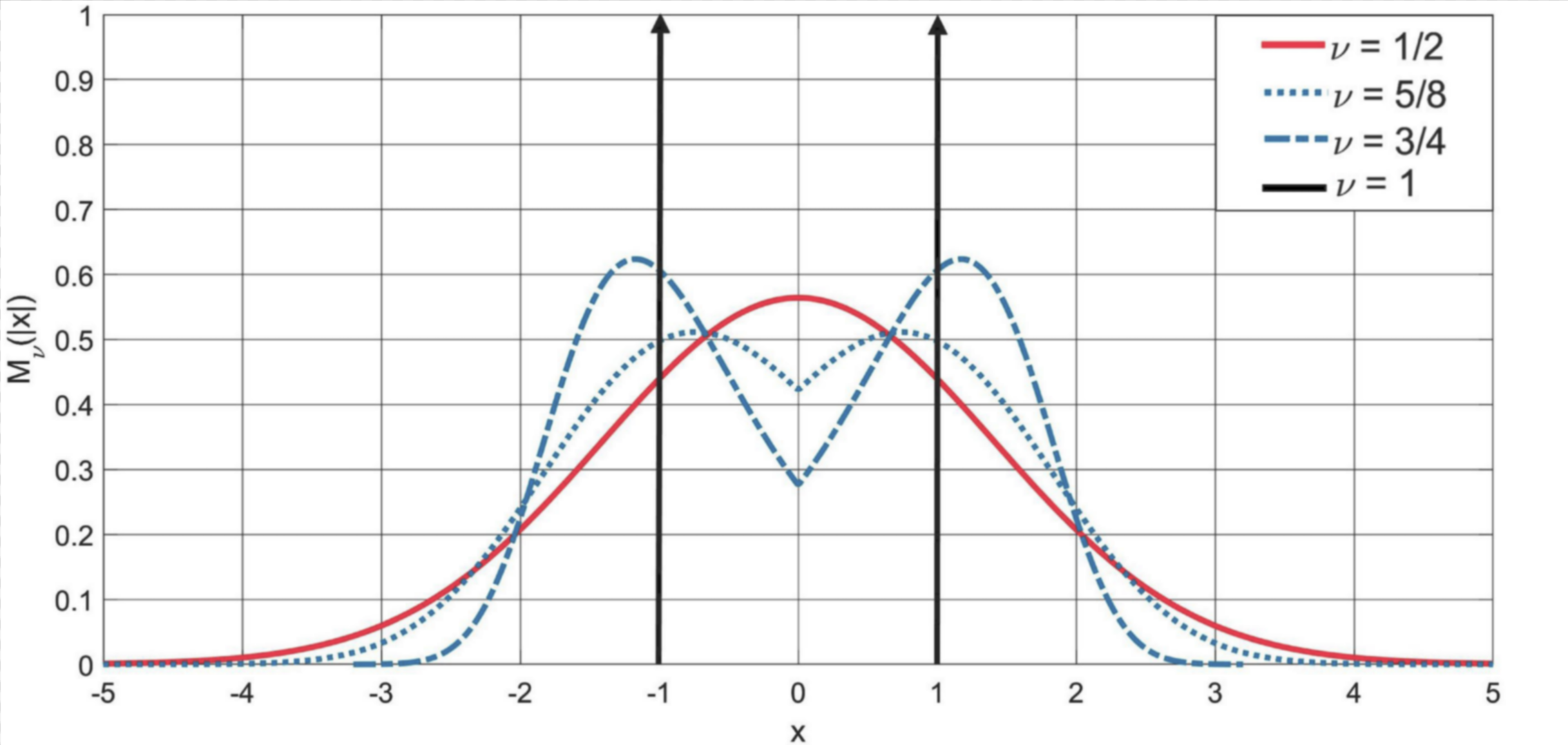}
%%	\caption
\vskip -0.45truecm
\caption{Plot of the symmetric $M-$Wright type function $M_{\nu}(|x|)|$
	 for $ 1/2 \leq \nu \leq 1$. Note that the $M$Wright function becomes a
a sum of two delta functions centered in $x=\pm 1$ for $\nu = 1$.}	 
\label{fig7}
\end{figure}
\newpage
\vsp
The  readers are invited to look the  YouTube  video 
by my former student Armando Consiglio whose title is ``Simulation of the  $M-$Wright function'', in which 
the author  has  shown  the evolution of this  function as the parameter $\nu$ changes between 0 and 0.85 in the  interval ($-5<x<+5$) of $\RR$ centered in $x=0$ represented 
 at fixed time $t=1$.
 %\vsp
Furthermore,  the interested readers can  have more details on the classical Wright functions
by  consulting  the recent survey by Luchko \cite{Luchko HFCA19} and references therein. 
\newpage 
\vsp
%% \subsection{the Wright $ \MM$-function in Two Variables.}
In view of time-fractional diffusion processes related to time-fractional diffusion equations 
it is worthwhile to introduce the 
function  in two variables
\begin{equation} 
  \MM_\nu(x,t):= t^{-\nu}\, M_\nu(xt^{-\nu})\,,\q 0<\nu < 1\,,\q x,t \in \RR^+ \,,
\label{F.51}
  \end{equation}   %% (F.51) 
  which defines a  spatial probability density in $x$ evolving in 
  time $t$ with self-similarity exponent $H=\nu$.
  Of course for $x\in \RR$ 
  we have to consider the symmetric version of the $M$-Wright function.
  obtained from 
  (\ref{F.51}) multiplying by $1/2$ and replacing $x$ by $|x|$.
 
%  \vsp 
   Hereafter we provide
  a list of the main properties of this function, 
  which can be derived from the Laplace and Fourier transforms 
  for the corresponding Wright $M$-function
  in one variable  presented in papers by Mainardi and recalled in the Appendix F
  of Reference \cite{Mainardi BOOK10}.

% \vsp
   For the Laplace transform
    of $\MM_\nu(x,t)$ with respect to $t>0$ and $x>0$  we get respectively:   
   \begin{equation}
   \L\left\{\MM_\nu (x,t);  t \to s \right\}
:= \int_0^\infty \e^{-st}\, t^{-\nu}\, M_\nu(x\, t^{-\nu})\, dt    
   = s^{\nu-1}\, \e^{\ds \, -xs^\nu}\,;
\label{F.52}
\end{equation}  
   %\eqno(F.52)$$
   %%%%%%%%%%%%%%%%%%%
\begin{equation}	
	\L\left\{\MM_\nu(x,t);  x \to s \right\}
:= \int_0^\infty \e^{-sx}\, t^{-\nu}\, M_\nu(x\, t^{-\nu})\, dx	
	=  E_{\nu, 1}\left( -s t^\nu \right)\,.
\label{F.53}
\end{equation}	
	%\eqno(F.53)$$
	%%%%%%%%%%%%%%%%%%%%%%%%%%%%%%%%%%%%%%%%%%%%%%%%%%%%%%%%%%%%%%%%%%%%%%%%%%%%%%%
  For the Fourier transforms with respect to the spatial variable $x$ we have
  for $\MM_\nu (x,t)$ with $x \in \RR^+$,
\bee
\begin{array}{ll}
 &\F_C\left\{\MM_\nu(x,t);  x \to \kappa \right\} :=
 \int_0^\infty \cos(\kappa x)\, t^{-\nu}\, M_\nu(x\, t^{-\nu})\, dx =
 E_{2\nu,1} (-\kappa^2 t^{2\nu})\,,
 \\ \\
 &\F_S\left\{\MM_\nu(x,t);  x \to \kappa \right\} :=
  \int_0^\infty \sin(\kappa x)\, t^{-\nu}\, M_\nu(x\, t^{-\nu})\, dx =
   \kappa^\nu\, E_{2\nu,\nu+ 1} (-\kappa^2 t^{2\nu})\,,
  \end{array}
  \label{F.54CS}
  \ee
  so that for the symmetric function 
$\MM_\nu(|x|,t)$  
  we get
  \bee
  \F\left\{\MM_\nu(|x|,t);  x \to \kappa \right\}=
2  \int_0^\infty \cos(\kappa x)\, t^{-\nu}\, M_\nu(x\, t^{-\nu})\, dx =
 2E_{2\nu, 1}\left( -\kappa^2 t^{2\nu} \right)\,.
\label{F.54}
\ee
%}

%%
Restricting our attention at the known analytic expressions of the $M_\nu$ functions versus $x$ at fixed time $t=1$ we recall the following results for some special rational values of the parameter $\nu$:
%\\
%$\nu= 1/4$
%\bee
%M_{1/4}(x) = \mathcal{E} (???)\,,
%\label{SZ-1}
%\ee
\\
$\nu=1/3$ (see Reference \cite{Mainardi BOOK10})
\bee
M_{1/3}(x) = 3^{2/3}\Ai(x/3^{1/3})\,,
\label{MC-Ai}
\ee
\\
$\nu =1/2$ (see Reference \cite{Mainardi BOOK10}
\bee
M_{1/2}(x) = \frac{1}{\sqrt{\pi}}\e^{-x^2 /4}\,,
\label{MC-GAUSS}
\ee
\\
$\nu =2/3$ (see Reference \cite{Hanyga PRSA02}) 
%  and Reference \cite{Liemert-Kleine JMP2015}) 
\bee
M_{2/3}(x) = 3^{-2/3}
\left[3^{1/3}\,x\, \Ai\left(x^2/3^{4/3}\right) -
3\Ai^\prime\left ( x^2/3^{4/3}\right)\right ]
\, \e^{-2x^3/27}\,.
\label{MC-Ai1}
\ee
%\\
%$\nu=3/4$
%\bee
%M_{3/4}(x) = \mathcal{E} (???)\,.
%\label{SZ-2}
%\ee
\vsn
In the above equations
%$\mathcal {E}$ and $\mathcal{E}^\prime$ denote the {\it Exponential integral} and its %%derivative,  and 
$\Ai$  and $\Ai^\prime$ denote the \emph{Airy function}
and its first derivative.

 %%%%%%%% 30 OCTOBER 2020 %%%

\section{Conclusions}
In this survey we have reviewed the main applications of the Mittag-Leffler function
in deterministic and stochastic processes. The analysis,
essentially based on contributions of the author by himself and co-authored with colleagues,
 is not intended to be exhaustive.
With all applications dealt in this survey,  we have tried to justify  the title of {The Queen Function of the Fractional Calculus}
 generally attached to the Mittag-Leffler function in the framework of researchers on fractional calculus. 

%%%%%%%%%%%%%%%%%%%%%%%%%%%%%%%%%%%%%%%%%%
\section*{Acknowledgments}
{The work of  the author
has been carried out in the framework of the activities of the National Group of Mathematical Physics (GNFM, INdAM).
The author would like to thank the anonymous reviewers for their helpful and constructive comments and Professor Bruce West for the invitation.}
%}
%In this section you can acknowledge any support given which is not covered by %the author contribution or funding sections. This may include administrative and %technical support, or donations in kind (e.g., materials used for experiments).}

%%%%%%%%%%%%%%%%%%%%%%%%%%%%%%%%%%%%%%%%%%
%\vvs \noindent
\subsection*{Funding}
% \funding
 {This research received no external funding.} 

\subsection*{Conflicts of interest}
{The author declares no conflict of interest} %Please add conflicts of interest.
%\appendixtitles{yes} %Leave argument "no" if all appendix headings stay EMPTY (then no dot is printed after "Appendix A"). If the appendix sections contain a heading then change the argument to "yes".
\%appendix
%\section{}

\section*{Appendix A: Author's  acquaintance with the Mittag-Leffler function since the late 1960's}
%%%%%%%%%%%%%%
I (the Author) was formerly acquainted with the Mittag-Leffler function 
from the pioneering 1947 paper by Gross on creep and relaxation in linear 
viscoelasticity. It was during my PhD sties at the University of Bologna under the supervision of Prof Caputo in the year 1969. Indeed  I was asked to apply 
in the framework of anelastic materials the
derivative of non-integer order introduced by Prof Caputo
in \cite{Caputo GJRAS67,Caputo BOOK69}. 
More recently this fractional derivative  was named after him thanks the suggestions of 
Gotrenflo and Mainardi \cite{Gorenflo-Mainardi CISM97} and 
Podlubny \cite{Podlubny BOOK99}.
I understood that the Mittag-Leffler function proposed by Gross
both in creep and relaxation processes could be used in the corresponding processes 
 in the fractional Zener model. Because Gross had computed and 
 plotted only the spectra, see Figure 1 in this article, I was interested to plot
 the Mittag-Leffler function on which I was addressed in the Third volume
 of the Handbook of the Bateman Project published in 1955
 \cite{Erdelyi BATEMAN}.
 Carrying out the plot of the Mittag-Leffler function
 $E_\alpha(-t^\alpha)$ using a Fortran program was not easy for me using its        
 power series representation, so I limited the time interval
 to $[0. 5]$ with ordinate in logarithmic scale.
 As far as I know this was the former plot of this function, see References
\cite{Caputo-Mainardi PAGEOPH1971,Caputo-Mainardi RNC1971} where 
the results of  my PhD thesis were published in 1971
 jointly with my supervisor.
 Later I was acquainted with the viscoelastic model by Rabotnov in 1948
 \cite{Rabotnov 48} and with the Russian school of Meskov and Rossikhin
    who used the so-called Rabotnov function, indeed related to the Mittag-Leffler function.
  and consequently with results similar to some extent to those by Caputo and Mainardi 
  published in 1971 \cite{Caputo-Mainardi PAGEOPH1971,Caputo-Mainardi RNC1971}.
    However, our work was totally independent from the Russian school
  (incidentally published in Russian),  
    as outlined in the Notes to the chapter 3 of my 2010 book, see pp.  74--76
    in \cite{Mainardi BOOK10}.
    More later, in the 1980  I was acquainted with the results by Bagley-Torvik and
    by Koeller that confirmed the relevant role of the Mittag-Leffler functions in linear viscoelastic models governed by constitutive laws of fractional order. 
    Once again their results crossed with those in
    References \cite{Caputo-Mainardi PAGEOPH1971,Caputo-Mainardi RNC1971}.
       However, I have to confess that, when in conferences of those years 
    I dealt with fractional derivatives in rheology, the audience remained 
    indifferent  if not hostile and laughable  so I left this topic preferring to 
    transfer my research interests to wave phenomena, in particular   on the effects of dissipation on linear dispersive waves.
    
%    \vsp
    Incidentally, in 1980's, I was also aware of the nice treatise by Harold T. Davis on the Theory of Linear operators published in 1936 \cite{Davis BOOK36}, where the author gave 
    information about the fractional calculus and the Mittag-Leffler function.
    It was my honor to publish a recent survey on the contributions by Davis and Gross
    (already recalled in Introduction), whom 
    I consider the pioneers of fractional relaxation processes
     in viscoelastic and dielectric materials \cite{Mainardi-Consiglio WSEAS2020}. 
     %%
%   \vsp
In the firsts years of 1990s under the  push of fractals, the relevance 
    of fractional derivatives (used not always in a correct way)  was outlined in several papers.  
For this I was induced to come back to fractional calculus.
It was just this occasion for me to devote my research interests  to the application
of fractional calculus in relaxation, oscillation phenomena governed by fractional ODEs
and diffusion, wave phenomena governed by fractional PDEs.
Once again I understood the relevance of the Mittag-Leffler functions
but also that of the Wright functions, both of them
classified as miscellaneous functions in the handbook of Bateeman project.
I must note that, as far as I known,  the Bateman  handbook was the only 
one published in English  to deal with these special function, 
and therefore accessible to me.

%\vsp 
The year 1994 was the golden year for me as far as my acquaintance with 
fractional calculus and related special functions is concerned.
Indeed I took profit by the acquaintance in three  different conferences
 with the late Prof Gorenflo and Prof. Nigmatullin (in Bordeaux, France),
 with Prof.  Podlubny and Prof. Caputo (in Astlanta, USA), and 
  with Prof Virginia Kiryakova and the late Prof. Stankovic (in Sofia, Bulgaria),
   among other authorities of the fractional calculus.
 But it was with Prof Gorenflo that I started a collaboration for more than 20 years
 (1995-2015) motivated by our common interest towards the potential of the Mittag-Leffler functions  in the applications of the  fractional calculus.

% \vsp
 Then, since 1997, I was interested in the emerging science of Econophysics  thanks mainly to my younger colleague Enrico Scalas. With Gorenflo,  Scalas and his student Raberto  we published some papers on the advent of Fractional Calculus in Econophysics, 
 see e.g. \cite{Mainardi PhysicaA00}   and  my historical survey in
 Mathematics~\cite{Mainardi MATHEMATICS20}.
 In 2007, on the occasion of the 80-birthday of Prof. Caputo, I published with Gorenflo a survey in Fractional Calculus and Applied Analysis \cite{Mainardi-Gorenflo FCAA07}
 where I took the liberty to propose for the Mittag-Leffler function
 the (successful) title of {\it the Queen Function of the Fractional Calculus}.
  Some years earlier, 
 Gorenflo had contacted the American Mathematical Society to give a specific number to the Mittag-Leffler function, that is 33E12, in the MSC classification.
 
% \vsp
 Gorenflo and I  promoted the Mittag-Lffler functions in several Conferences and 
 Workshops in all the world. In particular, I would like to recall
 my lectures in India (under invitation of Prof Mathai, director of the Center 
 of Mathematical Sciences, in Brazil (under invitation of Prof Edmundo Capelas de Oliveira,
 Campinas University) and in US (under invitation of Prof. Karniadakis, Brown University,
 see my course  \cite{Mainardi FC-BROWN-LN15}).    

%\vsp
I like to outline my gratitude to Professor Michele Caputo (1927) and Rudolf Gorenflo
(1930--2017) 
 for having provided me with useful advice in  earlier and later times,
 respectively. 
%%%
 It is my pleasure to enclose a  photo showing
 the author between them, taken in Bologna, April 2002. %%  (R.G. at my right, M.C. at my le 
 %\vskip -1truecm 
 \begin{figure}   %% [H]
\centering
\includegraphics[width=6.0truecm]{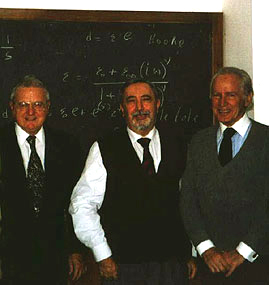}
\caption{F. Mainardi between R. Gorenflo (left) and M. Caputo (right).}
\end{figure}
% \begin{center}
%% \vskip -0.5truecm
%%%%%%%  \includegraphics[width=.90\textwidth]{gorenflo-mainardi-caputo.eps}
%%%%%%%\includegraphics[width=10.5truecm]{gorenflo-mainardi-caputo.eps}
% %% \includegraphics[width=10.5truecm]{fmws-figACK.jpg}
% \includegraphics[width=6.0truecm]{NEW_fmws-figACK.jpg}
%\end{center}
%\vskip -0.2truecm
%\centerline{F. Mainardi between R. Gorenflo (left) and M. Caputo (right).}
% %\centerline{Bologna, April 2002}
%\vsp
% Prof. Caputo  introduced me to the fractional calculus during my PhD thesis (1969--1971).
% Prof. Gorenflo has collaborated actively with me and with my students on  several papers for % 20 years   since 1995.
 Unfortunately,  I lost Gorenflo's  guidance and collaboration in 2015 when
 he suffered strong health troubles that led him to his death on 20 October 2017 at 87 years.  
 He was Emeritus Professor of Mathematics at the Free University of Berlin since his retirement in 1998.
 
 %%%%
 Nowadays  I am quite interested to promote the special functions of the Mittag-Leffler and Wright type with the second edition of the treatise by Gorenflo et al.
 \cite{GKMR BOOK20} and my surveys 
 \cite{Mainardi WSEAS2020,Mainardi-Consiglio MATHEMATICS20},
 including the present review.
 
 \section*{Appendix B:  Some notes on the numerical computation of the functions of the Mittag-Leffler type}
 
 In the Appendix A we have outlined the possible priority of the first plot of a
  Mittag-Leffler function in time domain in the 1971 paper by Caputo and Mainardi
  \cite{Caputo-Mainardi RNC1971}.
  
  Yu N. Rabotnov, in his works on viscoelastity, see 
 \eg \cite{Rabotnov 48,Rabotnov BOOK69,Rabotnov TABLES69}
 %% \cite{Rabotnov 48,Rabotnov 69,Rabotnov TABLES69}
 introduced the function of time $t$, depending on  two parameters, that he denoted by  $\alpha \in (-1,0]$
 (related to the type of viscoelasticity)
 and $\beta \in \RR$; 
 $$ R_\alpha(\beta, t) := 
 t^{\alpha}\, \sum_{n=0}^\infty \frac{\beta^n\, t^{n(\alpha +1)}}{\Gamma[(n+1)(\alpha+1)]}\,,\q t\ge 0\,.
 \eqno(B.1)
 $$
 Rabotnov,  almost surely unaware of the Mittag-Leffler function, 
  referred to this function  
 to as the {\it fractional exponential function} 
 noting that for $\alpha =0$ it reduces
 to the standard exponential $\exp (\beta t)$. 
  In the Russian literature such a function is mostly referred to as  the  {\it Rabotnov function}.  
 The relation of this function with  the Mittag-Leffler function in two parameters is thus obvious:   
 $$ R_{\alpha}(\beta, t) =  t^\alpha \, E_{\alpha+1,\alpha+1}\left(\beta \,t^{\alpha+1}\right )\,.\eqno(B.2)$$
 However, it is trivial  to show that
 $$R_{\alpha}(\beta, t) = \frac{1}{\beta}\,
 \dfrac{d}{dt} E_{\alpha +1}\left(\beta\,t^{\alpha +1} \right)\,. \eqno(B.3)$$
We have kept the notation of Rabotnov (except for the Russian letter denoting this function)  according to whom the main parameter $\alpha$  is
shifted with respect to ours,  so it appears confusing.
 However, in view of the fact that his 
constant $\beta$ is negative, we have a complete equivalence with our theory of fractional relaxation
discussed in  Section  {\bf 3}, namely with Eq. (8), if we set in (16)
 $(\alpha+1)\to \alpha $ and $\beta \to -1$.   
 The author, even when aware of some tables of this function published 
 by Rabotnov and his equipe in 1969
\cite{Rabotnov TABLES69}
 did never  use them: it could be interesting 
 to compare them with the numerical computations available nowadays.
 
 As far as we know,  the first paper entirely devoted to numerical analysis of the Mittag-Leffler function was that by  Gorenflo, Loutchko and Luchko
published in 2002, see {Gorenflo-Loutchko-Luchko FCAA02}
The authors gave details for setting a program with MATHEMATICA,
that then was translated in into a MATLAB code by Podlubny
\cite{Podlubny MATLAB06}, freely available at the   WEB Site of  MATLAB Central (File exchange).

 %% Computation of the Mittag-Leffler function and its derivatives. 
 %% Fract. Calc. Appl. Anal. 5(2002), 491-518

Also Hilfer with Seybold  was interested in producing routines
 for computing the Mittag-Leffler function
and its inverse 
\cite{Seybold-Hilfer FCAA05,Hilfer-Seybold ITSF06,Seybold-Hilfer SIAM08}  
but unfortunately we get only the procedures in  their  articles
but not any routine even if promised in C e MATLAB  in \cite{Seybold-Hilfer SIAM08}.

A routine in FORTRAN90 (mlfv.f90)
 was carried out by
Davide Verotta  for the evaluation of the two parameter Mittag-Leffler function
\cite{Verotta 10a,Verotta 10b} 
and then used and improved  by Mainardi and Spada in their 2011 paper on fractional viscoelastic models \cite{Mainardi-Spada EPJ-ST11}.

Diethelm in his nice book on Fractional Calculus 
\cite{Diethelm BOOK10} has given details
 for the computation of the Mittag-Leffler function using rational Pad\'e approximants.

We mention  the  routine in MATLAB proposed  in 2008 by Yang-Quan Chen 
\cite{Chen MATLAB08}
and the note of 2014  by Val\'erio and Machado \cite{Valerio-Machado CNSNS14}.

As far as we know the most efficient and accurate routine 
for the computation of the Mittag-Lefler function in the whole complex plane
(that I use with my collaborators) is that by Garrappa
\cite{Garrappa MATLAB14},  
illustrated in an issue of SIAM Journal of Numerical Analysis 
\cite{Garrappa SIAM15}, and based on the numerical inversion of the Laplace transform.
This routine allows us  to compute the Mittag-Leffler function with up to 3 independent parameters 
(the so-called Prabhakar function) that reads
$$  E_{\alpha, \beta}^\gamma (z) =
\sum_{n=0}^{\infty} \frac{(\gamma)_n}{n!\Gamma(\alpha n+ \beta)}\, z^n\,, \quad 
\Re \alpha > 0,\, \beta \in \CC, \, \gamma>0, \eqno(B.4)$$
where
$  (\gamma)_n = \gamma (\gamma+1)\dots (\gamma + n -1)= 
{\Gamma(\gamma+n)}/{\Gamma (\gamma)}$ denotes the Pochhammer symbol.
%%%
The Prabhakar function is fundamental in dealing with dielectric relaxation processes
governed by the Havriliak-Negami law, see \cite{Mainardi-Garrappa JCP15},
\cite{Garrappa-Mainardi-Maione FCAA16} and the more recent survey
by Giusti et al. \cite{Giusti FCAA20}.

Garrappa with his collaborators has also considered the numerical evaluation of 
matrix Mittag-Leffler functions (which however exits from the purpose of this survey) in a paper \cite{GarrappaPopolizio2018} which also discusses the numerical computations of derivatives of the scalar Mittag-Leffler function. 

The MATLAB routine by Garrappa has been translated  into Fortran and R languages by others, as stated in  Garrappa's home page 
 {\tt https://www.dm.uniba.it/members/garrappa}
at the Department of Mathematics, University of Bari, Italy.

% REFERENCES BIBLIOGRAPHY %%%
%\reftitle
%{\hl{References--refs. 4,5,11,41-46,51,54,57,73 are not cited}} %a lot of references %(4,5,11,41-46,51,54,57,73) are not cited in the main text, please cite them in the main text. %Please also make sure the citation and the reference are cited in numberic order. 

\end{document}